\documentclass{amsart}
\usepackage{graphicx}
\usepackage{amsmath}
\usepackage{amssymb}
\usepackage{bm}
\usepackage{hyperref}
\theoremstyle{plain}
\begingroup
 
\newtheorem{theorem}{Theorem} 
\newtheorem{assumption}{Assumption} 

\newtheorem{corollary}{Corollary}
\newtheorem{lemma}{Lemma}
\endgroup
\newcommand{\nwc}{\newcommand}
\nwc{\qref}[1]{(\ref{#1})}
\nwc{\cadlag}{c\`{a}dl\`{a}g}
\nwc{\la}{\label}
\nwc{\nn}{\nonumber}
\nwc{\ee}{\text{e}}
\nwc{\Z}{\mathbb{Z}}
\nwc{\C}{\mathbb{C}}
\nwc{\E}{\mathbb{E}}
\nwc{\R}{\mathbb{R}}
\nwc{\N}{\mathbb{N}}
\nwc{\Mn}{\mathbb{M}_N}
\nwc{\PP}{\mathcal{P}}
\nwc{\sub}{\mathcal{S}}
\nwc{\M}{\mathcal{M}}
\nwc{\law}{\stackrel{\mathcal{L}}{\rightarrow}}
\nwc{\eqd}{\stackrel{\mathcal{L}}{=}}
\nwc{\vp}{\varphi}
\nwc{\Vp}{\Phi}
\nwc{\ve}{\varepsilon}
\nwc{\veps}{\varepsilon}
\nwc{\eps}{\ve}
\nwc{\dnto}{\downarrow}
\nwc{\nsup}{^{(n)}}
\nwc{\ksup}{^{(k)}}
\nwc{\jsup}{^{(j)}}
\nwc{\nksup}{^{(n_k)}}
\nwc{\inv}{^{-1}}
\nwc{\argmin}{\mathrm{argmin}}
\nwc{\argmax}{\mathrm{argmax}}
\nwc{\Stoch}{\mathfrak{Q}}
\nwc{\ad}{\mathrm{ad}}
\nwc{\Ad}{\mathrm{Ad}}
\nwc{\diag}{\mathrm{diag}}
\nwc{\dd}{\mathfrak{d}}
\nwc{\gln}{\mathfrak{g}}
\nwc{\unn}{\mathfrak{u}}
\nwc{\lotri}{\mathfrak{l}}
\nwc{\GLM}{GL(m,\C)}
\nwc{\proj}{P}
\nwc{\Tr}{\mathrm{Tr}}
\nwc{\Id}{\mathrm{Id}}
\nwc{\son}{\mathfrak{so}(n)}
\nwc{\so}{\mathfrak{so}}
\nwc{\SO}{\mathsf{SO}}
\nwc{\SOn}{\mathsf{SO}(n)}
\nwc{\Symm}{\mathrm{Symm}}
\nwc{\SymmN}{\mathsf{Symm}(N)}
\nwc{\MM}{\mathcal{M}}
\nwc{\TT}{\mathbb{T}}
\nwc{\Rq}{\mathbb{R}^q}
\nwc{\GL}{\mathsf{GL}}
\nwc{\schwartz}{\mathcal{S}}
\nwc{\laplacian}{\triangle}
\nwc{\eigen}{\psi}
\nwc{\one}{\mathbf{1}}
\nwc{\id}{I}
\nwc{\Lap}{\mathcal{L}}
\nwc{\Sz}{Sz\'{e}kelyhidi}
\nwc{\Poincare}{Poincar\'{e}}
\nwc{\Mobius}{M\"{o}bius}
\nwc{\Arnold}{Arnol'{d}}
\nwc{\Muller}{M\"{u}ller}
\nwc{\Sverak}{\v{S}ver\'{a}k}
\nwc{\Berard}{B\'{e}rard}
\nwc{\Ampere}{Amp\`{e}re}
\nwc{\Scal}{\mathrm{Scal}}
\nwc{\Ric}{\mathrm{Ric}}
\nwc{\nnn}{\mathbf{n}}
\nwc{\bb}{\mathbf{b}}
\nwc{\Holder}{H\"{o}lder}
\nwc{\Levy}{L\'{e}vy}
\nwc{\T}{\mathbb{T}}
\nwc{\Varphi}{\Phi}
\nwc{\Ito}{It\^{o}}
\nwc{\vol}{\mathrm{vol}}
\nwc{\Vol}{\mathrm{Vol}}
\nwc{\rhoeq}{\rho_{\mathrm{eq}}}
\nwc{\matn}{\mathbb{M}}
\nwc{\grad}{\mathrm{grad}}
\nwc{\Lie}{\mathcal{L}}
\nwc{\Hess}{\mathrm{Hess}}
\nwc{\res}{\varepsilon}
\nwc{\ressq}{\tilde{\res}}
\nwc{\HH}{\mathcal{H}}
\nwc{\hern}{\mathbb{H}(n)}
\nwc{\On}{\mathbb{O}(n)}
\nwc{\DD}{\mathcal{D}}
\newcommand{\RN}[1]{%
  \textup{\uppercase\expandafter{\romannumeral#1}}%
}

\renewcommand{\Re}{\mathop{\rm Re}\nolimits}
\renewcommand{\Im}{\mathop{\rm Im}\nolimits}

\theoremstyle{definition}

\theoremstyle{remark}
 
\numberwithin{equation}{section}
\numberwithin{figure}{section}
\begin{document}
\title{Local rigidity of Siegel Brownian motion}
\author{Tianmin Yu}
\address{Department of Mathematics, Northwestern University, 2033 Sheridan Rd., Evanston, IL 60201. }

\email{tianmin.yu@northwestern.edu}
%
%

%\thanks{This work has been supported in part by the National Science Foundation (DMS 2107205). GM also acknowledges partial support from JSPS KAKENHI grants JP18H01124 and JP20K20884 and expresses his gratitude to the Research Institute for Mathematical Sciences at Kyoto University and the Department of Mathematics, Kyushu University during the completion of this work.}

%\subjclass[2010]{58G32,60D05,30.0X}
%\keywords{Brownian motion on Riemannian manifolds, Dyson Brownian motion, Symplectic geometry, Symmetric space}

\begin{abstract}
    We study the large $N$ behavior of Siegel Brownian motion for $\beta\geq 1$,
started from a fixed finitely supported empirical measure with rational
weights. Uniformly for $(\log N)^{-2}\leq t\leq \frac13\log N$, we show
that the particle counts of half-lines differ from their deterministic
predictions by $O(\log N)$ with high probability. This gives simultaneous
rigidity of order $(\log N)/N$ in regions where the limiting density is
bounded below. At a fixed positive time with a regular rightmost edge,
the largest particle lies at most $O(N^{-2/3}(\log N)^{2/3})$ below the
edge and $O(N^{-2/3}(\log N)^2)$ above it with high probability.
The proof uses a stochastic Hamilton--Jacobi equation for the normalized
logarithmic potential. Random characteristics control the counting error,
while deterministic characteristics and a comparison at two imaginary
heights yield the exterior edge estimate.
\end{abstract}

\maketitle
%\tableofcontents
\section{Introduction}
In a recent work~\cite{menon2025siegel}, Menon and the author considered Brownian motion on the Siegel upper half-space $\mathcal H_N$ and the evolution of orbits for a stabilizer subgroup of $\mathrm{Sp}(N,\mathbb R)$. This resulted in an autonomous interacting particle system with a hyperbolic interaction term. With proper scaling, the system can be described by the following SDE (with a different $\beta$ convention from the original one):
\begin{align}\label{sde:siegel}
    \mathrm d\sigma_i(t)=&\frac1N\big(\coth\sigma_i(t)+\sum_{j\neq i}\frac{\sinh\sigma_i(t)}{\cosh\sigma_i(t)-\cosh\sigma_j(t)}\big)\mathrm dt\nonumber\\
    &+\sqrt{\frac2{N\beta}}\mathrm dB_i(t)\,,&&i=1,...,N.
\end{align}
This SDE serves as a hyperbolic version of Dyson Brownian motion, and is named Siegel Brownian motion. This process $\{\sigma_i(t)\}_{i=1}^N$ also turned out to be a special type of radial Heckman--Opdam process introduced  by Schapira in~\cite{schapira2007heckman}. Recent work by Huang and McSwiggen~\cite{huang2024asymptotics} investigated the large deviation of general radial Heckman--Opdam process. Mundo~\cite{mundo2026geometric} investigated the geometric construction of this type of processes.

The large-$N$ limit of this system is of fundamental interest in random matrix theories. Dyson Brownian motion, unitary Brownian motion, and Siegel Brownian motion arise in Euclidean, compact, and noncompact settings, respectively, and the drift terms are all derivatives of log volume of isospectral submanifolds, resulting in $\frac1x,\cot x,\coth x$ type interactions. 
C\'epa and L\'epingle~\cite{cepa2001brownian} studied related particle systems and discussed the hyperbolic interaction.

Local rigidity of particles is a fundamental problem in the analysis of  this types of interacting particle systems, and is also the first step in the three-step strategy proposed by Erd\H{o}s and Yau~\cite{erdHos2017dynamical}. In this framework, the resolvent is the main object to be analyzed. From a dynamical point of view, the resolvent usually satisfies a type of stochastic Burgers' equation, which can be analyzed through method of characteristics.
In the works~\cite{adhikari2020dyson,adhikari2023local,huang2019rigidity}, this method was used to prove several rigidity estimates for general Dyson Brownian motion models and unitary Brownian motions.

The main contribution of this paper is to establish the local law and rigidity estimate up to poly-logarithmic factors for Siegel Brownian motion. In particular, the uniform bulk rigidity bound is of order $\frac{\log N}N$ which matches the order known for equilibrium log-gases~\cite{bourgade2022optimal}. At the edge, we proved the rigidity of order $O(N^{-\frac23}(\log N)^2)$(outside) and $O(N^{-\frac23}(\log N)^\frac23)$(inside), which has the correct power order in $N$.

The method used in this paper is the study of the smooth counting function, which is the primitive of the resolvent. The smooth counting function satisfies a stochastic Hamilton--Jacobi equation, which can be analyzed through a stochastic characteristic method. Bounds on the smooth counting function directly provide the estimates of the type $\gamma_{i-C\log N}<\sigma_i<\gamma_{i+C\log N}$, where $\sigma_i$ and $\gamma_i$ are the position and classical position of the $i$-th particle, yielding local rigidity estimates of the type $|\sigma_i-\gamma_i|<\frac {C\log N}{N\rho(\gamma_i)}$ in the bulk.
Outside the limiting support, a sharper estimate bounds the distance of the largest particle from the right edge.

The novelty of this method mainly comes from integrable structures in the Hamilton--Jacobi equation that provide much finer controls compared with existing methods analyzing Burgers' equation. The regularity of the smooth counting function also avoids several divergence issue that might appear in resolvents.

The proofs in this paper are of a PDE style and mainly uses stochastic analysis to obtain estimates. We first prove an a priori bound on the smooth counting function, which leads to rigidity in the bulk and inside the edge. Then we use the a priori bound to obtain a finer control and prove the rigidity outside of the edge.

\subsection{Main results}
In this paper, we use $C$ to denote large constants independent of $N$, and its meaning may vary line by line, and $N$ is supposed to be sufficiently large.

We consider an initial empirical measure of the following form that is independent of $N$,
\begin{align}
    \mu_0(\sigma)=\frac1{\sum_{i=1}^rn_i}\sum_{i=1}^rn_i\delta_{a_i}(\sigma)\,,
\end{align}
for $\{n_i\}_{i=1}^r\subset \mathbb N$ and $\{a_i\}_{i=1}^r\subset [0,\infty)$. We use $M_\mu=\max_i a_i$ to denote the largest $a_i$.

For $\mu_N(\cdot,0)=\frac1N\sum_{i=1}^N\delta_{\sigma_i(0)}$ to be the same as $\mu_0$, we need $N$ to be a multiple of $(\sum_{i=1}^rn_i)/\mathrm{g.c.d}(n_1,...,n_r)$.

We consider $\beta\geq 1$ which is the correct threshold introduced in~\cite{menon2025siegel} that ensures the SDE~\eqref{sde:siegel} to have a strong solution and the particles are not colliding. We also introduce $T_0=(\log N)^{-2}$ and $T_1=\frac13\log N$ to provide the range of time to be considered.

Let $\mathcal S^{++}=\{z=E+\bm i \eta:E>0,\eta\in(0,\pi)\}$ be the double positive strip,
we define the complex valued function $S_N(z,t)=\frac1N\sum_{i=1}^N\log(\cosh z-\cosh\sigma_i(t))$ for $z\in\mathcal S^{++}$. $S_N(z,t)$ satisfies a stochastic Hamilton--Jacobi equation
\begin{align}\label{SDE:S_N}
    &\mathrm dS_N(z,t)+\frac12[\partial_zS_N(z,t)]^2\mathrm dt=\frac{N+1}{2N}\mathrm dt+\frac1{N}R_N(z,t)\mathrm dt +\frac1N\mathrm dL_N(z,t)\\
    &\text{ where }R_N(z,t)=(\frac1\beta-\frac12)\partial^2_zS_N(z,t)-\frac12\coth(z)\partial_zS_N(z,t)\,,\\
    &\text {and }\mathrm dL_N(z,t)=\sqrt{\frac2{\beta N}}\sum_{i=1}^N\frac{-\sinh\sigma_i(t)}{\cosh z-\cosh\sigma_i(t)}\mathrm dB^i_t\,.
\end{align}
We also define the resolvent $m_N(z,t)=\partial_zS_N(z,t)=\frac1N\sum_{i=1}^N\frac{\sinh z}{\cosh z-\cosh \sigma_i(t)}$. It is a hyperbolic type Stieltjes transform of the empirical measure.

Let $S_\infty(z,t)$ be the solution (defined in~\ref{lem:m_inf-sol}) to the limiting equation
\begin{align}\label{eq:HJ-S_inf}
    &\partial_tS_\infty(z,t)+\frac12(\partial_zS_\infty(z,t))^2=\frac12\,,
\end{align}
with initial value $S_\infty(z,0)=S_N(z,0)$. We also define $m_\infty(z,t)=\partial_zS_\infty(z,t)$, and $\Delta_N(z,t)=S_N(z,t)-S_\infty(z,t)$. $S_\infty(z,t)$ and $m_\infty(z,t)$ are known to be explicitly solvable and unique through method of characteristics.

Once the limit $S_\infty(z,t)$ is known, 
consider the (double) counting functions $\mathfrak n_N(E,t),\mathfrak n_\infty(E,t)$ for $E\geq 0$, defined as
\begin{align}
    \mathfrak n_N(E,t)=&\frac1N\sum_{i=1}^N\Big(\bm 1_{\{\sigma_i(t)\leq E\}}+\bm 1_{\{-\sigma_i(t)\leq E\}}\Big)=1+\frac1N\sum_{i=1}^N\bm 1_{\{\sigma_i(t)\leq E\}}\,,\\
    \mathfrak n_\infty(E,t)
        =&2-\frac1\pi\lim_{\eta\to 0^+}\Im(S_\infty(E+\bm i\eta,t))\,.
\end{align}
The limit of the smooth counting function $S_N(E+\bm i\eta,t)$ as $\eta\to 0^+$, is known to characterize the counting function $\mathfrak n_N(E,t)$.

The main theorems we have are the following:
\begin{theorem}\label{thm:counting}
    For any $D>19$, there exists a constant $C$ depending on $D,\beta$, such that
    \begin{align}
        \sup_{E\geq 0}|\mathfrak n_N(E,t)-\mathfrak n_\infty(E,t)|<C\frac{ \log N+t^{-\frac12}}N\,,\text{ uniformly in } t\in(0,T_1]\,,
    \end{align}
    holds with probability at least $1-N^{-D}$. In particular, if $t>T_0=\frac1{(\log N)^2}$ the bound is just $C\frac{\log N}N$.    
\end{theorem}

This theorem comes from the following control on the function $\Delta_N(z,t)$.
\begin{theorem}\label{thm:bound-Delta}
    For $D>19$, there exists a constant $C$ depending on $D,\beta$ only, such that the following event happens with probability at most $N^{-D}$:
    \begin{align}
    \sup_{z\in(0,N^2)\times(N^{-2},\pi)\atop t\in[0,T_1]}
    |\Delta_N(z,t)|\geq C\frac{\log N}N\,.
\end{align}
\end{theorem}

The following corollary is a direct result of Theorem~\ref{thm:bound-Delta} and Cauchy integral formula for derivatives.
\begin{corollary}\label{coro:resolvent-convergence}
With probability at least $1-N^{-D}$, we have
\begin{align}
    &\sup_{t\in[0,T_1]}
    |m_N(z,t)-m_\infty(z,t)|< 2C\frac{\log N}{N\eta}\,,\nonumber\\
    &\text{ uniformly in }z\in(0,N]\times(2N^{-2},\frac\pi2)\,,
\end{align}
where $\eta=\Im(z)$. 
\end{corollary}

Theorem~\ref{thm:counting} leads to the rigidity of particle indices. Let $\gamma_i(t)$ be the classical position of $\mathfrak n_\infty(E,t)$, defined by
\begin{align}
    &\gamma_i(t)=\inf\{E\geq0:\mathfrak n_\infty(E,t)\geq1+\frac iN\}
    \,,&&i\in(0,N]\,,\\
    \text{ and }&\gamma_i(t)=0\,,&&i\leq 0\,,\\
    &\gamma_i(t)=\infty\,,&&i>N\,.
\end{align}
Then we have
\begin{corollary}\label{coro:index-rigidity}
    With probability at least $1-N^{-D}$, 
    \begin{align}
        \gamma_{i-C\log N}(t)< \sigma_i(t)< \gamma_{i+C\log N}(t)\,,
    \end{align}
    holds uniformly for all $t\in[T_0,T_1]$, $i=1,...,N$.
\end{corollary}

If a point $\gamma_i(t)$ is in the bulk, this implies the particle local rigidity 
\begin{align}
    |\sigma_i(t)-\gamma_i(t)|< C\frac{\log N}{N\rho}\,,
\end{align}
where $\rho$ is any number comparable with the density in between $\gamma_{i-C\log N}(t)$ and $\gamma_{i+C\log N}(t)$.

At the edge, Corollary~\ref{coro:index-rigidity} leads to rigidity inside the edge. Consider a right soft edge $E_0$ at time $t_0$ with a square root density $\rho(E,t_0)\sim \sqrt{(E_0-E)^+}$, then we have a lower bound on $\sigma_N(t)$:
\begin{align}
    E_0-\sigma_N(t_0)<E_0-\gamma_{N-C\log N}(t_0)<\tilde C(\frac{\log N}N)^{\frac23}\,.
\end{align}

To control the particle outside of the edge, we need a more detailed control on the imaginary part of $S_N(z,t)$. Bounds of the order $\frac{\log N}{N}$ still allows $\log N$ number of particles to go out, and a bound of the order less than $\frac1N$ is needed to prove that $\sigma_N(t)-E_0$ is not too large. We have the following theorem:
\begin{theorem}\label{thm:edge-estimate}
    With probability at least $1-N^{-D}$, 
    \begin{align}
        \sigma_N(t_0)-E_0<CN^{-\frac23}(\log N)^{2}\,,
    \end{align}
    for some constant $C$ independent of $N$.
\end{theorem}

This theorem comes from the following estimate.
\begin{theorem}\label{thm:edge}
    For a given $z=E_0+\kappa+\bm i\eta$,
    with $\kappa\in(0,r)$ and $\eta\in(4N^{-2},c\kappa)$ for some $r,c>0$ independent of $N$, the following holds with probability at least $1-N^{-D}$,
    \begin{align}
        |\Im(\Delta_N(z,t_0))|<\frac CN(\frac{\sqrt{\log N}\eta}{\kappa}+\frac{(\log N)^2\sqrt{\kappa}}{N\eta^2}+\sqrt{\frac{(\log N)^2\sqrt{\kappa}}{N\eta^2}})\,.
    \end{align}
\end{theorem}

With a scaling $\kappa=a N^{-\frac23}(\log N)^2$ and $\eta=a^\frac12 N^{-\frac23}(\log N)^\frac32$ for $a$ independent of $N$ but large enough, the bound of Theorem~\ref{thm:edge} is not exceeding $\frac{0.03}N$. This excludes any particle on the right of $E_0+2\kappa$.

As an example,
for $\mu_0=\delta_0$, the initial resolvent is
$m_\infty(z,0)=\coth(z/2)$. For $t\ll 1$, the density $-\frac1\pi\Im(m_\infty(E+\bm i0^+,t)$ 
is approximately semi-circle with radius $\sqrt{8t}$. For $t\gg 1$, the density is approximately uniform on $[0,t]$. The rightmost edge is given by
\begin{align}
    E_0(t)=2\,\mathrm{arcsinh}(\sqrt{\frac t2})+\sqrt{t(t+2)}\,,
\end{align}
and the density vanishes as a square root at $E_0(t)$.

\subsection{Organization of this paper}
In section~\ref{sec:pre} we review some basic results about method of characteristics, in particular those with a random driving field, which will help the path-wise analysis in building up the a priori bound. We also introduce a Ward-type estimate for this hyperbolic type of resolvent.

In section~\ref{sec:apriori} we establish the a priori bound and prove Theorem~\ref{thm:bound-Delta}, Theorem~\ref{thm:counting}. The main subtlety comes from the inversion of characteristics. 
For a fixed terminal point $z$, the inverse characteristic has a random starting point, so the resulting stochastic integrand need not be adapted. 
To deal with this issue we consider a grid $\{\xi_i\}$ of initial points such that at least one of $\xi_i$ is close to the initial point $\xi$ of $z$. We are able to provide a union bound on the values of $\Delta_N(z_i,t)$ through stochastic analysis, then with some regularity condition checked, we can prove that $z$ is sufficiently close of $z_i$ so $\Delta_N(z,t)$ is also controlled.

In section~\ref{sec:edge} we prove Theorem~\ref{thm:edge} and Theorem~\ref{thm:edge-estimate}. We trace $\Delta_N(z,t)$ along the deterministic characteristic of $m_\infty(z,t)$. The a priori bound provide a finer control of $\Im(\Delta_N(z,t))$.
Near the right edge, we use the a priori bound to control the remainder terms and estimate the quadratic variation of the imaginary martingale part.

Proofs of some Lemmas are in the Appendix~\ref{sec:appendix}.

\subsection{Acknowledgements}
The author thanks Antonio Auffinger for helpful discussions related to this problem. 
The author used ChatGPT (OpenAI, GPT-6 Astra Pro) for mathematical discussions, for help with exposition, and to draft the proofs of 
Lemmas~\ref{lem:nevanlinna-boundary-estimates} and ~\ref{lem:two-H-func}. The author takes full responsibility for the final statement and proof.

\section{Notation and preliminaries}\label{sec:pre}
For $\mu$ a probability measure on $\mathbb R_{\geq 0}$, we define the smooth counting function and the resolvent as:
\begin{align}
    S_\mu(z)=&\int_{\mathbb R} \log(\cosh z-\cosh x)\mathrm d\mu(x)\,,\\
    m_\mu(z)=&\partial_zS_\mu(z)=\int_{\mathbb R} \frac{\sinh z}{\cosh z-\cosh x}\mathrm d\mu(x)\,,
\end{align}
where we take the branch of $\log $ such that for $z\in\mathbb C^+$, $\Im(\log z)\in(0,\pi)$.
We mainly work with $z\in\mathcal S^{++}$, while $S_\mu$ and $m_\mu$ should be regarded as functions defined on $\mathcal S^+=\{E+\bm i\eta:E\in\mathbb R,\eta\in(0,\pi)\}$, obtained by $S_\mu(-\bar z)=\overline{S_\mu(z)}$.

Special for $\mu_N(t)=\frac1N\sum_{i=1}^N\delta_{\sigma_i(t)}$, we use $S_N(z,t)$ and $m_N(z,t)$ to denote its smooth counting function and resolvent.

The limit $S_\infty(z,t)$ is defined as the solution to the deterministic Hamilton--Jacobi equation~\ref{eq:HJ-S_inf}, its derivative therefore is a solution to the complex Burgers' equation
\begin{align}
    \partial_tm_\infty(z,t)+m_\infty(z,t)\partial_zm_\infty(z,t)=0\,,
\end{align}
where the initial value $S_\infty(z,0)$ and $m_\infty(z,0)$ are the same as $S_N(z,0)$ and $m_N(z,0)$.

\subsection{Method of characteristics}
Associated with any probability measure $\mu_t(\cdot)$ on $\mathbb R_{\geq 0}$ weakly continuous in $t$, and its resolvent $m(z,t)$,
we can introduce a characteristic $Z(\xi,t)$ being the solution to the ODE for $\xi\in S^{++}$, up to its first exit from $\partial\mathcal S^{++}$:
\begin{align}
    \partial_t Z(\xi,t)=m(Z(\xi,t),t)\,,&&Z(\xi,0)=\xi\,.
\end{align}
\begin{lemma}\label{lem:char-inverse-exist}
    For any $\mu_t$ weakly continuous in time, $Z(\xi,t)$ has a well-defined inverse map $\Xi:(z,t)\to\xi$ for all $t>0$ and $z\in\mathcal S^{++}$ that takes values in $\mathcal S^{++}$. 
    
    Furthermore, $\Xi(z,t)$ is holomorphic in $z$, and $\Im(\Xi(z,t))>\Im(z)$.
\end{lemma}
The proof is in Appendix~\ref{pf:char-inverse-exist}.

Tracing $m_\infty(z,t)$ along the characteristic of itself, $m_\infty(Z_\infty(\xi,t),t)$ is a constant in time and is equal to $m_\infty(\xi,0)$, so the characteristics become straight lines in $t$. 
\begin{lemma}\label{lem:m_inf-sol}
For any initial probability measure $\mu(\cdot)$ on $\mathbb R$ with compact support, equation~\ref{eq:HJ-S_inf} has a unique solution with initial condition $S_\infty(z,0)=\int\log(\cosh z-\cosh\sigma)\mathrm d\mu(\sigma)$ given by:
\begin{align}
    S_\infty(z,t)=S_\infty(\xi,0)+\frac t2m_\infty(\xi,0)^2+\frac t2\,,
\end{align}
where $\xi=\Xi(z,t)$ is the unique solution in $\mathcal S^{++}$ to the equation $z=\xi+tm_\infty(\xi,0)$.

Furthermore, the resolvent $m_\infty(z,t)$ and characteristic $Z_\infty(\xi,t)$ are 
\begin{align}
    Z_\infty(\xi,t)=&\xi+tm_\infty(\xi,0)\,,\qquad m_\infty(z,t)=m_\infty(\xi,0)\,.
\end{align}
\end{lemma}
The proof is in Appendix~\ref{pf:m_inf-sol}.

\subsection{Ward type estimate}
The Stieltjes transform $F(z)=\int\frac1{z-x}\mathrm d\mu(x) $ of a measure $\mu$ on $\mathbb R$ satisfies ward estimates:
\begin{align}
    |F^{(k)}(z)|\leq \frac{k!|\Im(F(z))|}{(\Im(z))^k}\,,&&k=1,2,...,
\end{align}
Similar estimates hold for the hyperbolic version $m_\mu(z)$.

For $z=E+\bm i\eta\in\mathcal S^{+}$, the following identity holds (notice that $\Im(\coth(\frac z2))<0$):
\begin{align}
    &|\coth(\frac z2)|^2=1-2\cot(\eta)\Im(\coth(\frac z2))\,,\\
    &|\partial_z\coth(\frac z2)|=-\frac{\Im(\coth(\frac z2))}{\sin\eta}\,.
\end{align}
Along with the fact that
\begin{align}
    \frac{\sinh z}{\cosh z-\cosh\sigma}=\frac12\coth(\frac{z+\sigma}2)+\frac12\coth(\frac{z-\sigma}2)\,,
\end{align}
we have the following Ward-type estimate.
\begin{lemma}\label{lem:control-mz}
    For any probability measure $\mu$ on $\mathbb R$, its resolvent $m_\mu(z)$ with $z\in\mathcal S^{++}$ satisfies
    \begin{align}
        |m_\mu(z)|^2\leq &1-2\frac{\Im(m_\mu(z))}{\eta}\,,\\
        |\partial_zm_\mu(z)|\leq&1-\frac{\Im(m_\mu(z))}\eta\,.
    \end{align}

    The first inequality also implies a naive bound
    \begin{align}
        -\Im(m_\mu(z))<\frac{6}\eta\,.
    \end{align}
\end{lemma}

For the term $\coth(z)m_N(z)$ that appears in $R_N(z,t)$ we also have
\begin{lemma}\label{lem:control-cothzm}
    For any resolvent $m_\mu(z)$, any $z\in\mathcal S^{++}$,
    \begin{align}
        |\coth(z)m_\mu(z)|\leq&1-\frac{\Im(m_\mu(z))}\eta\,.
    \end{align}    
\end{lemma}

For the local martingale term $\mathrm dL_N(z,t)$, we have the following Lemmas for its quadratic variations.
\begin{lemma}\label{lem:control-martingale}
    For $z\in\mathcal S^{++}$, the (norm) quadratic variation of the complex-valued local martingale term $\mathrm dL_N(z,t)$ satisfies
    \begin{align}
        \langle\mathrm dL_N(z,t),\mathrm dL^*_N(z,t)\rangle\leq \frac2\beta(1-2\frac{\Im(m_N(z))}{\eta})\mathrm dt\,.
    \end{align}

    The quadratic variation of $\mathrm d\Im(L_N(z,t))$ satisfies
    \begin{align}
        \langle\Im( \mathrm dL_N(z,t))\rangle=&\frac1\beta \Big(\Re\big[\partial_zm_N(z,t)+\coth(z)m_N(z,t)\big]\nonumber\\
        &-\frac{\Im(\sinh(z)m_N(z,t))}{\Im(\cosh(z))}\Big)\mathrm dt\,,
    \end{align}
    and there exists a universal constant $C$ such that
    \begin{align}
        \langle\Im( \mathrm dL_N(z,t))\rangle\leq&\frac C\beta \big(1-\frac{\Im(m_N(z))}{\eta}\big)\mathrm dt\,.
    \end{align}
    
\end{lemma}

\section{An a priori bound}\label{sec:apriori}
In this section we provide a bound on $\Delta_N(z,t)$, which satisfies a transport equation
\begin{align}
    &\mathrm d\Delta_N(z,t)+\bar m(z,t)\partial_z\Delta_N(z,t)\mathrm dt
    =\big(\frac{1}{2N}+\frac1NR_N(z,t)\big)\mathrm dt+\frac1N\mathrm dL_N(z,t)\,,
\end{align}
where 
\begin{align}
    \bar m(z,t)=\frac12(m_\infty(z,t)+m_N(z,t))\,.
\end{align}

Let $\bar Z(\xi,t)$ be the characteristic of $\bar m(z,t)$. For any $\xi\in\mathcal S^{++}$, consider $\Delta_N(\bar Z(\xi,t),t)$ before $\bar Z(\xi,t)$ leaves $\mathcal S^{++}$, it satisfies that
\begin{align}\label{eq:Delta_N-integral}
    \Delta_N(\bar Z(\xi,t),t)=&\int_0^t\Big(\frac{1}{2N}\mathrm ds+\frac1NR_N(\bar Z(\xi,s),s)\mathrm ds+\frac1N\mathrm dL_N(\bar Z(\xi,s),s)\Big)\,.
\end{align}

\subsection{Stopping time bound}
Given $\xi$, we fix a threshold $N^{-3}$ and introduce a stopping time
\begin{align}
    \tau_\xi=\inf\{t>0:\Im(\bar Z(\xi,t))<N^{-3}\}\,.
\end{align}
%where $\tau_\xi=0$ in the case the set in the infimum is empty.

Consider the stopped process $\Delta_N(\bar Z(\xi,t\wedge\tau_\xi),t\wedge\tau_\xi)$, we have the bound\begin{lemma}\label{lem:control-R-L}
    There exists a constant $C$ depending on $\beta$ only, such that for any $\xi\in\mathcal S^{++},\Im(\xi)>N^{-3},t\in[0,T_1]$,
    \begin{align}
        &\Big|\int_0^{t\wedge\tau_\xi}R_N(\bar Z(\xi,s),s)\mathrm ds\Big|<C\log N\,,\\
        &\Big\langle\int_0^{t\wedge\tau_\xi}\mathrm dL^*_N(\bar Z(\xi,s),s),\int_0^{t\wedge\tau_\xi}\mathrm dL_N(\bar Z(\xi,s),s)\Big\rangle
        <C\log N\,.
    \end{align}
\end{lemma}
\begin{proof}
    By Lemma~\ref{lem:control-mz} and Lemma~\ref{lem:control-cothzm},
    \begin{align}
        \Big|\int_0^tR_N(\bar Z(\xi,s),s)\mathrm ds\Big|<&C\int_0^t\big(1-\frac{\Im(m_N(\bar Z(\xi,s),s))}{\Im(\bar Z(\xi,s))}\big)\mathrm ds\nonumber\\
        <&C\int_0^t\big(1-2\frac{\Im(\bar m(\bar Z(\xi,s),s))}{\Im(\bar Z(\xi,s))}\big)\mathrm ds\,.
    \end{align}

    For the process $\Im(\bar Z(\xi,t))$ before leaving $\mathcal S^{++}$, it satisfies that
    \begin{align}
        \frac{\mathrm d}{\mathrm dt}\Im(\bar Z(\xi,t))=\Im(\bar m(\bar Z(\xi,t),t))\,,
    \end{align}
    so we find
    \begin{align}
        -\int_0^{t}\frac{\Im(\bar m(\bar Z(\xi,s),s))}{\Im(\bar Z(\xi,s))}\mathrm ds=&-\log\frac{\Im(\bar Z(\xi,{t}))}{\Im(\bar Z(\xi,0))}<\log \pi-\log\Im(\bar Z(\xi,{t}))\,.
    \end{align}

    Similarly, the local martingale term has a quadratic variation controlled by Lemma~\ref{lem:control-martingale}
    \begin{align}
        \Big\langle\int_0^t\mathrm dL_N(\bar Z(\xi,s),s),\int_0^t\mathrm dL^*_N(\bar Z(\xi,s),s)\Big\rangle<&\frac2\beta\int_0^t(1-2\frac{\Im(m_N(\bar Z(\xi,s),s))}{\Im(\bar Z(\xi,s))})\mathrm ds\nonumber\\
        <\frac2\beta\int_0^t(1-4\frac{\Im(\bar m(\bar Z(\xi,s),s))}{\Im(\bar Z(\xi,s))})\mathrm ds
        <&\frac2\beta(t+\log \pi-4\log (\Im(\bar Z(\xi,t)))\,.
    \end{align}

    Replace $t$ by $t\wedge \tau_\xi$, since $\Im(\bar Z(\xi,t\wedge\tau_\xi))\geq N^{-3}$, and $t\leq T_1=\frac13\log N$, we conclude that each of the terms are bounded by $C\log N$.
\end{proof}        
Now the stopped local martingale $\int_0^{t\wedge\tau_\xi}\mathrm dL_N(\bar Z(\xi,s),s)$ is a martingale, so
Doob's maximal inequality along with Chernoff bound gives that
for any $M>1$, 
\begin{align}
    \mathbb P\big[\sup_{t\in[0,T_1]}|\int_0^{t\wedge\tau_\xi}\mathrm dL_N(\bar Z(\xi,s),s)|\geq 4MC\sqrt{\log N}\big]<\exp(-M^2)\,.
\end{align}
For any $D>0$, take $M=\sqrt{2D\log N}$, we have
\begin{corollary}\label{coro:grid-single-bound}
    Fix $\xi\in\mathcal S^{++}$ with $\Im(\xi)>N^{-3}$. For any $D>0$, there exists a constant $C$ depending on $D,\beta$, such that
    \begin{align}
        \mathbb P[\sup_{t\in[0,T_1]}|\Delta_N(\bar Z(\xi,t\wedge \tau_\xi),t\wedge\tau_\xi)|\geq C\frac{\log N}N]<N^{-2D}\,.
    \end{align}
\end{corollary}

\begin{proof}
    By equation~\eqref{eq:Delta_N-integral}, we already have the bounds on the two integrals.
\end{proof}

% \begin{proof}
%     By equation~\eqref{eq:Delta_N-integral}, we already have the bounds on the two integrals, while the initialization error has a direct bound by the following:
%     \begin{align}
%         \Delta_N(z,0)=&\int\log(\cosh z-\cosh\sigma)(\mathrm d\mu_N(\sigma,0)-\mathrm d\mu_\infty(\sigma,0))\nonumber\\
%         =&\int\log(\frac{\cosh z-\cosh\sigma}{\cosh z-1})(\mathrm d\mu_N(\sigma,0)-\mathrm d\mu_\infty(\sigma,0))\,.
%     \end{align}

%     For $z=E+\bm i\eta$, with $\eta\in (N^{-3},\pi)$, $|\sigma|\in[0,M_\mu]$, where $M_\mu$ is bound of the support of $\mu_N(\cdot,0)$ and $\mu_\infty(\cdot,0)$, the following holds
%     \begin{align}
%         \frac{\eta^2}{\pi^2}\leq |\frac{\cosh z-\cosh\sigma}{\cosh z-1}|\leq \exp(M_\mu)\frac{\pi^2}{\eta^2}\,,
%     \end{align}
%     implying that for some $C>0$,
%     \begin{align}
%         |\log(\frac{\cosh z-\cosh\sigma}{\cosh z-1})|<C \log N\,,
%     \end{align}
%     so
%     \begin{align}
%         |\Delta_N(z,0)|<C\frac{\log N}N\,.
%     \end{align}

%     Combine this with Lemma~\ref{lem:control-R-L} and the probabilistic bound, we proved the claim.
% \end{proof}
%up to the time $\bar Z(\xi,t)$ hits $\partial\mathcal S^{++}$.
\subsection{Inversion of characteristics}
To pass the bound of $\Delta_N(\bar Z(\xi,t\wedge\tau_\xi),t\wedge\tau_\xi)$ to some arbitrary $z$'s, we use the following strategy.

We take a uniform grid $\{\xi_i\}_{i=1}^{N^{19}}$ in $(0,N^4]\times(N^{-3},\pi)$, so that for each $\xi\in (0,N^4)\times(N^{-3},\pi)$, there exists a $\xi_i$ in the grid such that $|\xi-\xi_i|<10 N^{-7.5}<N^{-7}$.

Use Corollary~\ref{coro:grid-single-bound} for each grid point, we obtain a union bound:
\begin{corollary}\label{coro:grid-union-bound}
    For any $D>19$, we have the union bound for some constant $C$ depending on $D,\beta$,
    \begin{align}
       &\mathbb P\Big[\bigcup_{\xi_i\in\mathrm{grid}}\big\{ \sup_{t\in[0,T_1],t\leq\tau_{\xi_i}}|\Delta_N(\bar Z(\xi_i,t),t)|>C\frac{\log N}{N}\big\}\Big]<N^{19-2D}<N^{-D}\,,
    \end{align}
\end{corollary}

Now for each $z\in(0,N^2]\times(N^{-2},\pi)$, $T\in[0,T_1]$, consider its inverse $\xi=\bar \Xi(z,T)$ which is random in $\mathcal S^{++}$. By Lemma~\ref{lem:control-mz}
\begin{align}
    \Im(\bar m(z,t))<\frac{6}{\Im(z)}\,,
\end{align}
and that $\Im(\bar Z(\xi,t))$ is decreasing in time, we have
\begin{align}
    \Re(\xi)<\Re(z)+T \frac{6}{\Im(z)}<N^2+6TN^2<N^3\,,
\end{align}
so $\xi$ is still in the range of the grid, and we can take a $\xi_i\in\mathrm{grid}$ such that $|\xi-\xi_i|<N^{-7}$.

\begin{lemma}\label{lem:char-derivative-bound}
    For any measure $\mu_t(\sigma)$ continuous in time, the characteristic $Z(\xi,t)$ of its resolvent $m_\mu(z,t)$ satisfies that
    \begin{align}
        |\frac{\partial Z(\xi,t)}{\partial\xi}|\leq\frac{\pi\exp(t)}{\Im(Z(\xi,t))}\,.
    \end{align}
\end{lemma}

\begin{proof}
    For the ODE
    \begin{align}
        \partial_tZ(\xi,t)=m_\mu(Z(\xi,t),t)\,,
    \end{align}
    by Lemma~\ref{lem:control-mz} we have
    \begin{align}
        \Big|\frac{\partial Z(\xi,t)}{\partial\xi}\Big|=&\Big|\exp\big(\int_0^t\partial_zm_\mu(Z(\xi,s),s)\mathrm ds\big)\Big|\nonumber\\
        \leq &\exp\big(\int_0^t|\partial_zm_\mu(Z(\xi,s),s)|\mathrm ds\big)\nonumber\\
        \leq &\exp\big(\int_0^t(1-\frac{\Im(m_\mu(Z(\xi,s),s)}{\Im(Z(\xi,s)})\mathrm ds\big)\,.
    \end{align}

    Since $\frac{\mathrm d}{\mathrm ds}\Im(Z(\xi,s))=\Im(m_\mu(Z(\xi,s),s)$, the integral can be integrated and we have
    \begin{align}
        \Big|\frac{\partial Z(\xi,t)}{\partial\xi}\Big|\leq\exp(t)\frac{\pi}{\Im(Z(\xi,t))}\,.
    \end{align}
\end{proof}
\begin{corollary}\label{coro:near-z-exist}
     For any $T\in[0,T_1]$, any measure $\mu_t(\sigma)$ continuous in time, any $z_1\in\mathcal S^{++}$ with $\Im(z_1)>N^{-2}$, let $\xi_1=\Xi(z_1,T)$ be its initial point. For any $\xi_2\in\mathcal S^{++}$ with $|\xi_1-\xi_2|<N^{-7}$, the characteristic $\{Z(\xi_2,t)\}_{t\in[0,T]}\in\mathcal S^{++}$, and it satisfies
     \begin{align}
         &\inf_{t\in[0,T]}\Im(Z(\xi_2,t))>\frac12\Im(z_1)\,,\\
         \text{ and }&\sup_{t\in[0,T]}|Z(\xi_1,t)-Z(\xi_2,t)|<N^{-4.5}\,.
     \end{align}

\end{corollary}
\begin{proof}
    We prove by contradiction. For $s\in[1,2]$, consider $\xi_s=(2-s)\xi_1+(s-1)\xi_2$.
    By continuity of the characteristics in time, if there is a time $t\in[0,T]$ such that $\Im(Z(\xi_s,t))\leq\frac12 \Im(z_1)$ for some $s\in[1,2]$, then
    there is a time $t_0\in[0,T]$ and $s_0\in[1,2]$ such that $\Im(Z(\xi_{s_0},t_0))=\frac12 \Im(z_1)$, and $\Im(Z(\xi_{s},t_0))\geq\frac12 \Im(z_1)$ for all $s\in[1,2]$.
    
    Since $\Im(Z(\xi_1,t))$ as a function of $t$ is decreasing, we have
    \begin{align}
        \Im(Z(\xi_1,t_0)-Z(\xi_{s_0},t_0))>\frac12\Im(z_1)>\frac12N^{-2}\,.
    \end{align}

    On the other hand, by Lemma~\ref{lem:char-derivative-bound},
    \begin{align}
        |Z(\xi_{s_0},t_0)-Z(\xi_1,t_0)|\leq&\int_1^{s_0}|\frac{\partial Z(\xi_s,t_0)}{\partial\xi}|\times|\xi_2-\xi_1|\mathrm ds\nonumber\\
        \leq &\frac{\pi\exp(t_0)}{\frac12\Im(z_1)}N^{-7}<N^{-4.5}\,,
    \end{align}
    which is a contradiction.

    So we have that for all $t\in[0,T]$, $s\in[1,2]$, $\Im(Z(\xi_s,t))>\frac12 \Im(z_1)$. In particular, this holds for $s=2$. Furthermore,
    \begin{align}
        |Z(\xi_1,t)-Z(\xi_2,t)|\leq&\int_1^{2}|\frac{\partial Z(\xi_s,t)}{\partial\xi}|\times|\xi_2-\xi_1|\mathrm ds\nonumber\\
        \leq &\frac{\pi\exp(T)}{\frac12\Im(z_1)}N^{-7}<N^{-4.5}\,.
    \end{align}
\end{proof}

\subsection{Proof of theorem~\ref{thm:bound-Delta}}\label{sec:apriori-pf1}
\begin{proof}[Proof of theorem~\ref{thm:bound-Delta}]
    For any $z\in(0,N^2]\times(N^{-2},\pi)$, and $T\in[0,T_1]$, 
    we can find its inverse $\xi=\bar \Xi(z,T)$ and a neighboring $\xi_i$ in the grid, such that $|\xi-\xi_i|<N^{-7}$.
    Corollary~\ref{coro:near-z-exist} guarantees that $\{\bar Z(\xi_i,t)\}_{t\in[0,T]}\in\mathcal S^{++}$, and $\bar Z(\xi_i,t)$ satisfies that $\inf_{t\in[0,T]}\Im(\bar Z(\xi_i,t))>\frac12 N^{-2}$ and $\sup_{t\in[0,T]}|\bar Z(\xi_i,t)-\bar Z(\xi,t)|<N^{-4.5}$.

    Let $z_i=\bar Z(\xi_i,T)$, by the naive bound $|m_N(z,t)|,|m_\infty(z,t)|<\frac6{\Im(z)}$,
    \begin{align}
        &|\Delta_N(z,T)-\Delta_N(z_i,T)|=|\int_0^{1}(m_N-m_\infty)((1-\alpha)z+\alpha z_i,T)(z-z_i)\mathrm d\alpha|\nonumber\\
        &\leq \frac{12}{\frac12\Im(z)}|z-z_i|<24N^{2}N^{-4.5}<N^{-2}\,.
    \end{align}

    Take $C$ to be the same as the one in Corollary~\ref{coro:grid-union-bound}, the event $|\Delta_N(z,T)|>2C\frac{\log N}N$ for any $z\in(0,N^2]\times(N^{-2},\pi)$ and $T\in[0,T_1]$ implies that there exists $\xi_i$ in the grid such that
    \begin{align}
        |\Delta_N(z_i,T)|=|\Delta_N(\bar Z(\xi_i,T\wedge \tau_{\xi_i}),T\wedge \tau_{\xi_i})|>C\frac{\log N}N\,,
    \end{align}
    which implies the union event.

    So we conclude that
    \begin{align}
        \mathbb P\Big[\sup_{z\in(0,N^2)\times(N^{-2},\pi)\atop t\in[0,T_1]}|\Delta_N(z,t)|>C\frac{\log N}N\Big]<N^{-D}\,,
    \end{align}
    for $D>19$ and some constant $C$ depending on $D,\beta$.
    
\end{proof}

\subsection{Smooth counting function}\label{sec:apriori-smooth}
In this section we prove the Theorem~\ref{thm:counting} and Corollary~\ref{coro:index-rigidity}.

Corollary~\ref{coro:index-rigidity} is a direct result of Theorem~\ref{thm:counting}.
\begin{proof}[Proof of Corollary~\ref{coro:index-rigidity}]\label{pf:coro:index-rigidity}
    Given that
    \begin{align}
        \sup_{E\geq 0,t\in[T_0,T_1]}|\mathfrak n_N(E,t)-\mathfrak n_\infty(E,t)|<C\frac{\log N}N\,,
    \end{align}
    we only need to evaluate $\mathfrak n_\infty(\sigma_i(t),t)$.
    \begin{align}
        \mathfrak n_\infty(\sigma_i(t),t)<\mathfrak n_N(\sigma_i(t),t)+C\frac{\log N}N=1+\frac iN+C\frac{\log N}N
    \end{align}
    implies that $\sigma_i(t)<\gamma_{i+C\log N}(t)$, because $\mathfrak n_\infty(E,t)$ is a non-decreasing function in $E$. Similar for the other side.
\end{proof}

Now we consider the approximation of counting function through smooth counting function. For given $E\in\mathbb R$, the idea is to perturb $E$ to a point $E+h+\bm i\eta\in\mathbb C^+$, and compare $\mathfrak n_N(E,t)$ with $\Im(S_N(E+h+\bm i\eta,t))$. By the monotonicity of the counting function $\mathfrak n_N(E,t)$ in $E$ and smoothness of $S_N(z,t)$ in $z$ we can approximate $\mathfrak n_N(E,t)$ in terms of $\Im(S_N(E+h+\bm i\eta,t))$.

We introduce the following Lemmas about counting functions:
\begin{lemma}\label{lem:counting-approx}   
    Let $g_0(z)=1-\frac1\pi\Im(\log z)$, for $z\in\mathbb C^+$, in which case 
    $g_0(E+\bm i\eta)=\frac12+\frac1\pi\arctan(\frac{E}{\eta})$. For any $h>0,\eta>0$,
    \begin{align}
        g_0(E-h+\bm i\eta)-\frac{\eta}{\pi h}\leq \bm 1_{\{E\geq 0\}}\leq g_0(E+h+\bm i\eta)+\frac{\eta}{\pi h}\,.
    \end{align}    
\end{lemma}
\begin{proof}
    We only prove the right part of the inequality. The left part can be proved similarly.
    \begin{enumerate}
        \item If $E<0$, $\bm 1_{\{E\geq 0\}}=0$, while $g_0(\cdot)\geq 0$ and $\frac\eta{\pi h}>0$.
        \item If $E\geq 0$, $\bm 1_{\{E\geq 0\}}=1$, but
        \begin{align}
            g_0(E+h+\bm i\eta)+\frac{\eta}{\pi h}=&\frac12+\frac1\pi(\arctan(\frac{E+h}{\eta})+\frac{\eta}h)\nonumber\\
            \geq &\frac12+\frac1\pi(\arctan(\frac{h}{\eta})+\arctan(\frac{\eta}h))=1\,.
        \end{align}
    \end{enumerate}
\end{proof}

\begin{lemma}\label{lem:logch-log}
    For $\sigma\in\mathbb R$, $z=E+\bm i\eta\in\mathcal S^{+}$, let
    \begin{align}
        g_1(z,\sigma)=&2-\frac1\pi\Im(\log(\cosh z-\cosh \sigma))\,,\\
        h(z,\sigma)=&g_1(z,\sigma)-g_0(z-\sigma)-g_0(z+\sigma)\,,
    \end{align}
    where the convention of $\log(\cosh z-\cosh \sigma)$ such that 
    \begin{align}
        0\leq \Im(\log(\cosh z-\cosh \sigma))\leq 2\pi\,.
    \end{align}

    Then 
    there exists a universal constant $C$ such that
    \begin{align}
        |h(E+\bm i\eta,\sigma)|<C\eta\,.
    \end{align}
\end{lemma}
\begin{proof}
    By the convention of the argument of $\log (z)$, we can check that $\lim_{\eta\to 0^+}h(E+\bm i\eta,\sigma)=0$. 

    We aim to prove that $\tilde h(z,\sigma)=\log(\cosh z-\cosh \sigma)-\log(z-\sigma)-\log(z+\sigma)$ has bounded $z$-derivative. In fact,
    \begin{align}
        \partial_z\tilde h(z,\sigma)=\Big[\frac12\coth(\frac{z-\sigma}{2})-\frac1{z-\sigma}\Big]+\Big[\frac12\coth(\frac{z+\sigma}{2})-\frac1{z+\sigma}\Big]\,.
    \end{align}

    We only need to prove that 
    \begin{align}
        |\frac12\coth(\frac {z}2)-\frac1{z}|<C
    \end{align}
    for $|\Im(z)|\leq\pi$, which is true since this function is analytic everywhere in the region $\mathcal D=\{E+\bm i\eta:E\in[-1,1],\eta\in[-\pi,\pi]\}$, while for $|E|>1$ and $|\eta|\leq\pi$,
    \begin{align}
        &|\coth\frac {E+\bm i\eta}2|^2=\frac{\cosh E+\cos\eta}{\cosh E-\cos\eta}\leq \frac{\cosh E+1}{\cosh E-1}< 5\,,\\
        &|\frac1{E+\bm i\eta}|<1\,.
    \end{align}    
\end{proof}
\begin{corollary}\label{coro:logch-counting}
    For every $\sigma,E\in\mathbb R$, $h>0$, $\eta\in(0,\pi)$, it holds that
    \begin{align}
         \bm 1_{\{E\geq \sigma\}}+\bm1_{\{E\geq -\sigma\}}\leq& g_1(E+h+\bm i\eta,\sigma)+C(1+\frac1h)\eta\,,\\
         \bm 1_{\{E\geq \sigma\}}+\bm1_{\{E\geq -\sigma\}}\geq& 
         g_1(E-h+\bm i\eta,\sigma)-C(1+\frac1h)\eta\,.
    \end{align}
    where $C$ is a universal constant.
    
    In particular, this implies that
    \begin{align}
    &\mathfrak n_N(E,t)\leq 2-\frac1\pi\Im(S_N(E+h+\bm i\eta,t))+C(1+\frac1h)\eta\,\\
    \text{ and }&\mathfrak n_N(E,t)\geq 2-\frac1\pi\Im(S_N(E-h+\bm i\eta,t))-C(1+\frac1h)\eta\,.
    \end{align}
\end{corollary}

We also need a short time bound on $m_\infty(z,t)$ to bound the difference
\begin{align}
    S_\infty(E+h+\bm i\eta,t)-\lim_{\eta'\to 0^+}S_\infty(E+\bm i\eta',t)\,.
\end{align}
\begin{lemma}\label{lem:m_inf-short-time}
    For any $z\in\mathcal S^{++}$, $t>0$, 
    \begin{align}
        |m_\infty(z,t)|< \sqrt{1+\frac2t}\,.
    \end{align}
\end{lemma}
\begin{proof}
    Let $\xi=\Xi_\infty(z,t)$. The deterministic characteristic is
    \begin{align}
        z=\xi+tm(\xi,0)\,.
    \end{align}
    Take imaginary part, we have
    \begin{align}
        -\frac{\Im(m(\xi,0))}{\Im(\xi)}=\frac{\Im(\xi)-\Im(z)}{t\Im(\xi)}<\frac1t\,.
    \end{align}

    On the other hand, Lemma~\ref{lem:control-mz} gives that 
    \begin{align}
        |m(\xi,0)|^2\leq 1-2\frac{\Im(m(\xi,0))}{\Im(\xi)}<1+\frac2t\,,
    \end{align}
    so
    $|m_\infty(z,t)|=|m(\xi,0)|<\sqrt{1+\frac2t}$.
\end{proof}

Finally, the following Lemma restricts the support of $\mu_\infty(\cdot,t)$ in a finite time.

\begin{lemma}\label{lem:mu-inf-sup}
    Given that $\mu_\infty(\cdot,0)$ is supported on $[0,M_\mu]$, and $t\in[0,T_1]$, then 
    \begin{align}
        \mu_\infty([\frac N2,\infty),t)=0\,,
    \end{align}
    which implies $\mathfrak n_\infty(\frac N2,t)=0$.
\end{lemma}

\begin{proof}[Proof of Theorem~\ref{thm:counting}]
    First consider the case $E\in[0,N]$.
    We prove the $\leq$ side. The other side comes similarly.

    Take $h=N^{-1},\eta=2N^{-2}$, Corollary~\ref{coro:logch-counting} gives that
    \begin{align}
        \mathfrak n_N(E,t)<2-\frac1\pi\Im(S_N(E+h+\bm i\eta,t)+4C\frac1N\,.
    \end{align}

    By Theorem~\ref{thm:bound-Delta}, with probability at least $1-N^{-D}$,
    \begin{align}
        &|\Im(S_N(E+h+\bm i\eta,t)-S_\infty(E+h+\bm i\eta,t))|\nonumber\\
        &<|\Delta_N(E+h+\bm i\eta,t)|<C\frac{\log N}N\,,
    \end{align}

    For $S_\infty(z,t)$, Lemma~\ref{lem:m_inf-short-time} implies
    \begin{align}
        &|S_\infty(E+h+\bm i\eta,t)-\lim_{\eta'\to 0}S_\infty(E+\bm i\eta',t)|\nonumber\\
        <&\sqrt{1+\frac2t}\sqrt{h^2+\eta^2}<\sqrt{1+\frac2t}\frac2N<(2+4t^{-\frac12})\frac1N\,.
    \end{align}

    Combining these estimates we have
    \begin{align}
        \mathfrak n_N(E,t)<&2-\frac1\pi\Im(S_N(E+h+\bm i\eta,t)+4C\frac1N\nonumber\\
        \leq &\mathfrak n_\infty(E,t)+\frac1\pi|S_\infty(E+h+\bm i\eta,t)-\lim_{\eta'\to 0}S_\infty(E+\bm i\eta',t)|\nonumber\\
        &+\frac1\pi|S_\infty(E+h+\bm i\eta,t)-S_N(E+h+\bm i\eta,t)|+\frac CN\nonumber\\
        <&\mathfrak n_\infty(E,t)+C\frac{\log N+t^{-\frac12}}N.
    \end{align}

    Now we have the following holds with probability at least $1-N^{-D}$,
    \begin{align}
        \sup_{E\in[0,N]}|\mathfrak n_N(E,t)-\mathfrak n_\infty(E,t)|<C\frac{\log N+t^{-\frac12}}N\,,\text{ uniformly in }t\in[0,T_1]\,,
    \end{align}
    and we can extend $E\in[0,N]$ to $E\geq 0$.

    In fact, since $\mathfrak n_N(E,t)$ is non-decreasing in $E$, Lemma~\ref{lem:mu-inf-sup} implies that for any $E>N$, $t\in[0,T_1]$,
    \begin{align}
        \mathfrak n_N(E,t)\geq\mathfrak n_N(N,t)>2-C\frac{\log N+t^{-\frac12}}N=\mathfrak n_\infty(E,t)-C\frac{\log N+t^{-\frac12}}N\,.
    \end{align}
    The other side of the inequality is actually 
    \begin{align}
        \mathfrak n_N(E,t)\leq2=\mathfrak n_\infty(E,t)\,,
    \end{align}
    so the proof is completed.

\end{proof}

\section{Outside of the edge}\label{sec:edge}
In this section, we work in the event
\begin{align}
    E_\Delta=\{\omega:\sup_{z\in(0,N^2)\times(N^{-2},\pi)\atop t\in[0,T_1]}
    |\Delta_N(z,t)|< C\frac{\log N}N\}\,,
\end{align}
which has probability at least $1-N^{-D}$ by theorem~\ref{thm:bound-Delta}.

Recall that the deterministic characteristic $Z_\infty(\xi,t)$ of $m_\infty(z,t)$ is given by
\begin{align}
    &Z_\infty(\xi,t)=\xi+tm(\xi,0)\,,\\
    &m_\infty(Z_\infty(\xi,t),t)=m(\xi,0)\,.
\end{align}

Our final goal is to obtain a better control of the imaginary part of $\Delta_N(z,t)=S_N(z,t)-S_\infty(z,t)$. 
To make a good use of this deterministic characteristic, we trace $S_N(z,t)$ along it. Recall that $S_N(z,t)$ satisfies the stochastic Hamilton-Jacobi equation~\eqref{SDE:S_N}. 
\begin{lemma}\label{lem:SN-in-minf}
    For $z\in\mathcal S^{++}$, $t>0$, $\xi=\Xi_\infty(z,t)$, we use $Z_s$ to denote $Z_\infty(\xi,s)$ for short, then
    \begin{align}
        \Delta_N(z,t)=\frac{1}{2N}t+I_1+I_2+I_3+I_4\,,
    \end{align}
    where
    \begin{align}
        I_1=&\frac1N\int_0^tR_\infty(Z_s,s)\mathrm ds\,,&&
        I_2=\frac1N\int_0^tR_\Delta(Z_s,s)\mathrm ds\,,\nonumber\\
        I_3=&-\frac12\int_0^t\big(m_N(Z_s,s)-m_\infty(Z_s,s)\big)^2\mathrm ds\,,&&
        I_4=\frac1N\int_0^t\mathrm dL_N(Z_s,s)\,,
    \end{align}
    and
    \begin{align}
        R_\infty(z,t)=&(\frac1\beta-\frac12)\partial_zm_\infty(z,t)-\frac12\coth(z)m_\infty(z,t)\,,\\
        R_\Delta(z,t)=&(\frac1\beta-\frac12)\partial_z^2\Delta_N(z,t)-\frac12\coth(z)\partial_z\Delta_N(z,t)\,,\\
        \mathrm dL_N(z,t)=&-\sqrt{\frac2{\beta N}}\sum_{i=1}^N\frac{\sinh\sigma_i(t)}{\cosh z-\cosh\sigma_i(t)}\mathrm dB^i_t\,.
    \end{align}
\end{lemma}

For these four terms, we need to estimate them one by one.
The integral $I_1$ is a major term that is deterministic. The integral $I_2$ and $I_3$ are stochastic but can be bounded through the a priori bound and turns out to be minor. The quadratic variation of the imaginary part of martingale term $I_4$ turns out to be the summation of a major deterministic term with a stochastic minor term that can be bounded using the a priori bound. We evaluate these term in the following subsections.

\subsection{Integral of the major terms}
$I_1$ is completely computable (in the sense that it only depends on the starting and final point $\xi$ and $z$).
\begin{lemma}\label{lem:edge-I1}
    $I_1=\frac1N\Big((\frac1\beta-\frac12)\log(1+tm'(\xi,0))-\frac12\log(\frac{\sinh  z}{\sinh\xi})\Big)$.
\end{lemma}
\begin{proof}
    From the equation $z=\xi+tm(\xi,0)$, one can compute
    \begin{align}
        \partial_zm_\infty(z,t)=\partial_zm(\xi,0)=\frac{\partial_\xi(m(\xi,0))}{(\frac{\partial z}{\partial\xi})}=\frac{m'(\xi,0)}{1+tm'(\xi,0)}\,.
    \end{align}
    Therefore,
    \begin{align}
        \int_0^t\partial_zm_\infty(Z_s,s)\mathrm ds=\int_0^t\frac{m'(\xi,0)}{1+sm'(\xi,0)}\mathrm ds=\log(1+tm'(\xi,0))\,.
    \end{align}

    For the other term, there is still a similar integrable structure:
    \begin{align}
        &\int_0^t\coth(Z_s)m_\infty(Z_s,s)\mathrm ds=\int_0^t\coth(Z_s)m(\xi,0)\mathrm ds\nonumber\\
        =&\int_0^t\coth(Z_s)\mathrm dZ=\log(\frac{\sinh  z}{\sinh\xi})\,.
    \end{align}
\end{proof}

\subsection{Control of the minor terms}
The terms $I_2$ and $I_3$ all depend on the value of $\Delta_N(z,t)$, so with the apriori bound they should be small.

We introduce the following Lemma specific for $Z_\infty(\xi,t)$:
\begin{lemma}\label{lem:int-eta-2}
    For $z\in\mathcal S^{++}$ and $\xi=Z^{-1}_\infty(z,t)$, let $\eta_s=\Im(Z_\infty(\xi,s))$, then
    \begin{align}
        \int_0^t\frac{1}{\eta_s^2}\mathrm ds=\int_0^t\frac{1}{(\eta_0+s\Im(m_\infty(\xi,0)))^2}\mathrm ds=\frac{t}{\eta\eta_0}\,,
    \end{align}
    where $\eta=\Im(z),\eta_0=\Im(\xi)$.
\end{lemma}
\begin{corollary}\label{coro:I2-I3}
    For $\omega\in E_\Delta$, $z=E+\bm i\eta\in\mathcal S^{++}$ with $\eta\in(3N^{-2},\pi)$, $E\in(0,N^2)$, the following bound holds for any $t\in[0,T_1]$
    \begin{align}
        |I_2|+|I_3|<C\frac{(\log N)^2t}{N^2\eta\eta_0}\,.
    \end{align}
\end{corollary}
\begin{proof}
    Using the bound on $\Delta_N(z,t)$ and Cauchy integral formula, for any $z=E+\bm i\eta\in\mathcal S^{++}$ with $\eta\in(2N^{-2},\frac\pi2)$, take a circle centered at $z$ with radius $\frac\eta2$, we have
    \begin{align}
        &|\partial_z\Delta_N(z,t)|<C\frac{\log N}{N\frac\eta2}\,,\\
        &|\partial_z^2\Delta_N(z,t)|<2!C\frac{\log N}{N(\frac\eta2)^2}\,.
    \end{align}

    For $I_3$ part, 
    \begin{align}
        |I_3|<&\frac12\int_0^t|\partial_z\Delta_N(z,t)|^2\mathrm ds<2C^2\frac{(\log N)^2}{N^2}\int_0^t\frac{\mathrm ds}{\eta_s^2}
        =2C^2\frac{(\log N)^2t}{N^2\eta\eta_0}\,.
    \end{align}

    For the term $\coth(z)\partial_z\Delta_N(z,t)$ in $I_2$, use that $|\coth(z)|<1+\frac2{\Im(z)}$, we have a similar bound.
\end{proof}

\subsection{The martingale term}
Since we are integrating along a deterministic characteristic, the reverting $z\to\xi$ is deterministic, so $\int_0^t\mathrm dL_N(Z(\xi,s),s)$ is a martingale.

To make a detailed control, we should only consider the imaginary part of $I_4$ which will actually go into the counting function. The real part of $I_4$ could not be controlled as good as the imaginary part.

% \begin{lemma}\label{lem:im-quadratic}
%     The following identity holds for the quadratic variation of the imaginary part of $\mathrm dL_N(z,t)$:
%     \begin{align}\label{eq:ImdL-quadratic}
%         \langle\Im( \mathrm dL_N(z,t))\rangle=&\frac1\beta \Big(\Re\big[\partial_zm_N(z,t)+\coth(z)m_N(z,t)\big]\nonumber\\
%         &-\frac{\Im(\sinh(z)m_N(z,t))}{\Im(\cosh(z))}\Big)\mathrm dt\,.
%     \end{align}
% \end{lemma}
\begin{lemma}\label{lem:im-control}
    For $\omega\in E_\Delta$, $z=E+\bm i\eta\in\mathcal S^{++}$ with $\eta\in(3N^{-2},\pi)$, $E\in(0,N^2)$, there exists a constant $C$ depending on $D,\beta$,
    \begin{align}
        \Big|\Big\langle\Im\big(\int_0^t\mathrm dL_N(Z_s,s)\big)\Big\rangle-\frac1\beta\Theta\Big|<C\frac{(\log N)t}{N\eta\eta_0}\,,
    \end{align}
    where
    \begin{align}
        \Theta=\mathrm{Re}\Big(\log\big((1+tm'(\xi,0))\frac{\sinh z}{\sinh \xi}\frac{\Im(\cosh \xi)}{\Im(\cosh z)}\big)\Big)\,.
    \end{align}
\end{lemma}
The proof is in Appendix~\ref{pf:im-control}.

\subsection{The hyperbolic derivative}
Recall that
\begin{align}
    N I_1=&(\frac1\beta-\frac12)\log(1+tm'(\xi,0))-\frac12\log(\frac{\sinh z}{\sinh \xi})\,,\\
    \Theta=&\Re\Big(\log\big((1+tm'(\xi,0))\frac{\sinh z }{\sinh \xi}\frac{\Im(\cosh \xi)}{\Im(\cosh z)}\big)\Big)\,.
\end{align}

We rewrite them into the form
\begin{align}
    N I_1=&(\frac12-\frac1\beta)\log(\frac{\partial(\cosh\xi)}{\partial(\cosh z)})+\frac1{2\beta}\log(\frac{\cosh^2\xi-1}{\cosh^2z-1})\,,\\
    \Theta=&-\Re\Big(\log\big(\frac{\partial(\cosh\xi)}{\partial(\cosh z)}\frac{\Im(\cosh z)}{\Im(\cosh \xi)}\big)\Big)\,.
\end{align}

For $Z_\infty(\xi,t)$, its inverse $\Xi_\infty(z,t)$ is a univalent holomorphic map from $\mathcal S^{++}$ to itself.
Under a change of variable $w=\cosh z$, $\phi=\cosh\xi$, and notice that $\cosh(\mathcal S^{++})=\mathbb C^+$, for every fixed $t$, the map $F(w)=\cosh(\Xi_\infty(\mathrm{arccosh}(w),t))$ is a univalent holomorphic map from $\mathbb C^+$ to itself. The hyperbolic derivative of $F$ is defined as
\begin{align}
    F^h(w)=F'(w)\frac{\Im(w)}{\Im(F(w))}\,,
\end{align}
so we find that
\begin{align}
    \Im(NI_1)=&(\frac12-\frac1\beta)\Im(\log(F^h(w)))+\frac1{2\beta}\Im(\log(\frac{F(w)^2-1}{w^2-1}))\,,\\
    \Theta=&-\Re(\log (F^h(w)))\,.
\end{align}
In particular, $\Theta\geq0$ by Schwarz--Pick Lemma.

There are several classical results about this type of functions. We use the following results from theory of Nevanlinna functions and harmonic analysis to help edge estimate.
\begin{lemma}
\label{lem:nevanlinna-boundary-estimates}
Let
\begin{align}
    \mathbb H := \{z\in\mathbb C:\operatorname{Im}z>0\},
\end{align}
and let \(F:\mathbb H\to\mathbb H\) be a nonconstant holomorphic
function. Suppose that \(F\) extends holomorphically across
\(I=(a,b)\), where \(1\le a<b<\infty\), and that
\begin{align}
    F(s)>1,\qquad s\in I.
\end{align}
For \(x\in I\), set
\begin{align}
    d(x):=\min\{x-a,b-x\},
    \qquad
    F^h(x,y):=
    \frac{yF'(x+iy)}{\operatorname{Im}F(x+iy)}.
\end{align}
Then, for every \(x\in I\) and \(0<y<d(x)/2\), with \(z=x+iy\),
\begin{align}
    \left|
        \operatorname{Arg}
        \frac{F(z)^2-1}{z^2-1}
    \right|
    &< \frac{4y}{d(x)},
    \label{eq:nevanlinna-ratio-arg}
    \\
    \left|\operatorname{Arg}F^h(x,y)\right|
    &< \frac{2y}{d(x)},
    \label{eq:nevanlinna-hyperbolic-arg}
    \\
    0\le -\log\left|F^h(x,y)\right|
    &< \frac{8y^2}{3d(x)^2}.
    \label{eq:nevanlinna-hyperbolic-log}
\end{align}
Here \(\operatorname{Arg}\) denotes the principal argument.
\end{lemma}

% \begin{lemma}\label{lem:disk-map}
%     Let $\mathbb D=\{z\in\mathbb C:|z|<1\}$ and
% $\mathbb D^+=\{z\in\mathbb D:\Im z>0\}$.
% Suppose that $f:\mathbb D\to\mathbb D$ is holomorphic and univalent,
% that $f((-1,1))\subset\mathbb R$, and that
% $f(\mathbb D^+)\subset\mathbb D^+$. Define
% \[
%     f^h(z):=f'(z)\frac{\Im z}{\Im f(z)},
%     \qquad z\in\mathbb D^+.
% \]
% Choose the continuous branch of $\mathrm{Arg} f^h(iy)$ that tends to zero
% as $y\downarrow0$. Then
% \begin{align}
%     \bigl|\mathrm{Arg} f^h(iy)\bigr|
%     &\le 4\arctan y<4y,
%     && 0<y<\sqrt{2}-1,
%     \label{eq:fh-argument}\\
%     \bigl|\log |f^h(iy)|\bigr|
%     &\le \log\frac{1+y^2}{1-y^2}<8 y^2,
%     && 0<y<\frac12.
%     \label{eq:fh-log-modulus}
% \end{align}
% \end{lemma}
% \begin{lemma}
% \label{lem:pick-argument-interval}
% Let $\mathbb H=\{z\in\mathbb C:\Im z>0\}$, and let
% $F:\mathbb H\to\mathbb H$ be holomorphic. Suppose that $F$ extends
% holomorphically across an interval $I=(a,b)$, where
% $1\le a<b\le\infty$, and that $F(t)>1$ for every $t\in I$.
% Then, for $z=x+i\eta$ with $x\in I$ and $\eta>0$,
% \[
%     \left|\operatorname{Arg}\frac{F(z)^2-1}{z^2-1}\right|
%     \le
%     4\arctan\frac{\eta}{d(x)}<4\frac\eta{d(x)}\,,
% \]
% where $d(x)=\min\{x-a,b-x\}$.
% \end{lemma}
The proof is generated by GPT-6 Astra Pro and is in Appendix~\ref{pf:nevanlinna}.

\subsection{The edge rigidity}
Now we are able to prove the estimate around the edge. %Assume that $E(t)$ is an edge of $m_\infty(z,t)$, defined by the condition $\frac{\partial}{\partial\xi}Z_\infty(\xi,t)=0$ when ${Z_\infty(\xi,t)=E(t)}$.
For a fixed $O(1)$ time $t_0>0$, assuming a right edge $E_0>0$:
\begin{assumption}
    There exists a $r_1>0$ such that $\Xi_\infty(z,t_0)$ analytically extends to $\{z=E_0+r_1(1+\zeta):|\zeta|<1\}$, with $\Xi_\infty(\bar z,t_0)=\overline{\Xi_\infty(z,t_0)}$, and $\Xi_\infty(E,t_0)$ takes real values if $E\in(E_0,E_0+2r_1)$.
\end{assumption}
Notice that this implies that the density $-\frac1\pi\Im(m_\infty(E+\bm i 0^+,t_0))$ goes to $0$.

Consider the function $F(w)=\cosh(\Xi_\infty(\mathrm{arccosh}(w),t_0))$, there also exists $r_2>0$, such that $F(w)$ analytically extends to $\{w=\mathrm{arccosh}(E_0)+r_2(1+\zeta):|\zeta|<1\}$, with $F(\bar w)=\overline{F(w)}$, and $F(x)$ takes real values for $x\in(\mathrm{arccosh}(E_0),\mathrm{arccosh}(E_0)+2r_2)$.

\begin{corollary}\label{coro:hyperbolic-bound}
%    By Lemma~\ref{lem:disk-map} and Lemma~\ref{lem:pick-argument-interval}, 
    For $w=\cosh(E_0)+x+\bm iy$, with $x\in(0,r_2)$ and $0<y<0.1x$
    \begin{align}
        &|\mathrm{Arg}(F^h(w))|<2\frac{y}{x}\,,&|\log|F^h(w)||<\frac83\frac{y^2}{x^2}\,,\\
        \text{ and }& |\mathrm{Arg}\frac{F(w)^2-1}{w^2-1}|<4\frac yx\,.
    \end{align}

    Thus, for $z=E_0+\kappa+\bm i\eta$, there exists $r_3>0$ and $C,c>0$ depending on $E_0$, such that for $\kappa\in(0,r_3)$, $0<\eta<c\kappa$, 
    \begin{align}
        &|\mathrm{Arg}(F^h(\cosh z))|<C\frac{\eta}{\kappa}\,,&|\log|F^h(\cosh z)||<C\frac{\eta^2}{\kappa^2}\,,\\
        \text{ and }&|\mathrm{Arg}\frac{F(\cosh z)^2-1}{\cosh z^2-1}|<C\frac \eta\kappa\,.
    \end{align}
\end{corollary}
\begin{proof}
    Apply Lemma~\ref{lem:nevanlinna-boundary-estimates} to $f(z)=F(\mathrm{cosh}(E_0)+x(1+z))$.
\end{proof}

We also assume the edge is a soft one (square root density):
\begin{assumption}
    In a circle centered at $E_0$, 
    $m_\infty(z,t_0)=A(z)-B(z)\sqrt{z-E_0}$, where $A,B$ are analytic and take real values on the real axis, and $c=B(E_0)>0$.    
\end{assumption}

Under this assumption, we have
\begin{lemma}\label{lem:edge-assumption}
    There exists $r_4>0$, such that for $|\kappa+\bm i\eta|<r_4$, 
    \begin{align}
        |\Im(m_\infty(E_0+\kappa+\bm i\eta,t_0))|>\begin{cases}
            \frac c4\frac{\eta}{\sqrt{\kappa+\eta}}\,,&\kappa>0\,,\\
            \frac c4 \sqrt{|\kappa|+\eta}\,,&\kappa<0\,.
        \end{cases}
    \end{align}
\end{lemma}

\begin{proof}[Proof of Theorem~\ref{thm:edge}]\label{pf:thm-edge}
    By Lemma~\ref{lem:SN-in-minf}, we have 
    \begin{align}
        \Im(\Delta_N(z,t_0))=\Im(I_1)+\Im(I_2)+\Im(I_3)+\Im(I_4)\,.
    \end{align}

    Consider $z=E_0+\kappa+\bm i\eta$ with $\kappa\in(0,r)$ where $r=\frac12\min\{r_1,r_2,r_3,r_4\}$ and $\eta<c\kappa$,
    Corollary~\ref{coro:I2-I3} gives that
    \begin{align}
        |\Im(I_2)+\Im(I_3)|<C\frac{(\log N)^2t}{N^2\eta\eta_0}\,,
    \end{align}
    and Corollary~\ref{coro:hyperbolic-bound} gives that
    \begin{align}
        |\Im(I_1)|<\frac CN\frac{\eta}{\kappa}\,.
    \end{align}
    Combine Corollary~\ref{coro:hyperbolic-bound} with Lemma~\ref{lem:im-control}, with Doob's maximal inequality and Chernoff bound, we know that
    \begin{align}
        \mathbb P\Big[E_\Delta\cap\Big\{|\Im(I_4)|>\frac1N\sqrt{4CD\log N(\frac{\eta^2}{\kappa^2}+\frac{\log N t}{N\eta\eta_0})}\Big\}\Big]<N^{-D}\,.
    \end{align}

    Union bound with the event $E_\Delta$, we know that with probability at least $1-2N^{-D}$, 
    \begin{align}
        |\Im(\Delta_N(z,t_0))|<&C(\frac{\eta}{N\kappa}+\frac{(\log N)^2t}{N^2\eta\eta_0})+\frac1N\sqrt{4CD\log N(\frac{\eta^2}{\kappa^2}+\frac{\log N t}{N\eta\eta_0})}\nonumber\\
        <&C(\frac{\eta}{N\kappa}+\frac{(\log N)^2t}{N^2\eta\eta_0})+\frac1N\sqrt{4CD\log N\frac{\eta^2}{\kappa^2}}+\frac1N\sqrt{4CD\frac{(\log N)^2 t}{N\eta\eta_0}}\,.
    \end{align}

    Now for $\eta_0=\Im(\xi)$, since $\xi=z-tm_\infty(\xi,0)=z-tm_\infty(z,t_0)$, with $\Im(m_\infty(z,t_0))<0$ and Lemma~\ref{lem:edge-assumption}, and $\eta<c\kappa$, we have
    \begin{align}
        \eta_0=\Im(\xi)>t|\Im(m_\infty(z,t_0))|>Ct\frac{\eta}{\sqrt{\kappa+\eta}}>\tilde Ct\frac{\eta}{\sqrt{\kappa}}\,,
    \end{align}
    so
    \begin{align}
        \frac{t}{\eta\eta_0}<C\frac{\sqrt\kappa}{\eta^2}\,.
    \end{align}

    In summary, we have
    \begin{align}
        |\Im(\Delta_N(z,t_0))|<\frac CN\Big(\frac{\sqrt{\log N}\eta}{\kappa}+\frac{(\log N)^2\sqrt{\kappa}}{N\eta^2}+\sqrt{\frac{(\log N)^2\sqrt{\kappa}}{N\eta^2}}\Big)\,.
    \end{align}
\end{proof}

A direct result of this theorem is the $O(\frac1N)$ bound of $\Im(\Delta_N(z,t_0))$.
\begin{corollary}\label{coro:Delta_bound_outedge}
    For any $\epsilon>0$, there exists $\kappa_\epsilon\in(0,r),\eta_\epsilon\in(0,c\kappa_\epsilon)$, such that for $z_1=E_0+\kappa_\epsilon+\bm i\eta_\epsilon$, $z_2=E_0+\kappa_\epsilon+2\bm i\eta_\epsilon$, with probability at least $1-2N^{-D}$,
    \begin{align}
        \max\{|\Im(\Delta_N(z_1,t_0))|,|\Im(\Delta_N(z_2,t_0))|\}<\frac\epsilon N\,.
    \end{align}
\end{corollary}
\begin{proof}
    We consider the scaling $\kappa=a N^{-\frac23}(\log N)^2$ and $\eta=\sqrt a N^{-\frac23}(\log N)^\frac32$. Such $\kappa$ and $\eta$ satisfies $\kappa\in(0,r)$ and $\eta<c\kappa$ for any $O(1)$ choice of $a$. In this case, 
    \begin{align}
        &|\Im(\Delta_N(z_1,t_0))|<\frac CN\Big(\frac2{ a^\frac12}+\frac1{a^\frac14}\Big)\,,\nonumber\\
        &|\Im(\Delta_N(z_2,t_0))|<\frac CN\Big(\frac3{ a^\frac12}+\frac1{2a^\frac14}\Big)\,.
    \end{align}
    So it suffices to
    take $a=\max\{(\frac{4C}\epsilon)^2,(\frac{4C}\epsilon)^4\}$.
\end{proof}

To exclude the particle outside of $E_0+CN^{-\frac23}(\log N)^2$, we introduce the following Lemma.
\begin{lemma}\label{lem:two-H-func}
    For $E>0,\sigma\geq 0$, $\eta\in(0,0.1)$, consider the function
    \begin{align}
        H(E,\eta,\sigma)=&2\Im(\log(\cosh(E+\bm i\eta)-\cosh\sigma))\nonumber\\
        &-\Im(\log(\cosh(E+2\bm i\eta)-\cosh\sigma))\,.
    \end{align}

    The following properties hold:
    \begin{enumerate}
        \item $H$ is increasing in $\sigma$ for $\sigma\in[0,\infty)$.
        \item $0<H(E,\eta,\sigma)\leq8(\frac{\eta}{E-\sigma})^3$ if $\sigma\in[0,E)$.
        \item $\pi>H(E,\eta,\sigma)\geq\pi-4(\frac{\eta}{\sigma-E})^3$ if $\sigma\in(E,\infty)$.
    \end{enumerate}
\end{lemma}
The proof is generated by GPT-6 Astra Pro and is in Appendix~\ref{pf:two-H-func}.

\begin{proof}[Proof of Theorem~\ref{thm:edge-estimate}]\label{pf:thm-edge-estimate}
Consider the function 
\begin{align}
    D(E,\eta)=2\Im(\Delta_N(E+\bm i\eta,t_0))-\Im(\Delta_N(E+2\bm i\eta,t_0))\,.
\end{align}
Take $\epsilon=0.01$ in Corollary~\ref{coro:Delta_bound_outedge}, we obtain that
\begin{align}
    |D(E_0+\kappa_{0.01},\eta_{0.01})|<\frac{0.03}N\,.
\end{align}

On the other hand, 
\begin{align}
    D(E,\eta)=&\int_0^\infty H(E,\eta,\sigma)\mathrm d(\mu_N(\sigma,t_0)-\mu_\infty(\sigma,t_0))\,.
\end{align}
Using integral by part, 
\begin{align}
    &\int_0^{E_0} H(E,\eta,\sigma)\mathrm d(\mu_N(\sigma,t_0)-\mu_\infty(\sigma,t_0))\nonumber\\
    =&H(E,\eta,E_0)\delta\mathfrak n(E_0,t_0)-H(E,\eta,0)\delta\mathfrak n(0,t_0)\nonumber\\
    &-\int_0^{E_0}(\partial_\sigma H(E,\eta,\sigma))\delta\mathfrak n(\sigma,t_0)\mathrm d\sigma\,,
\end{align}
where $\delta\mathfrak n(E,t)=\mathfrak n_N(E,t)-\mathfrak n_\infty(E,t)$.

Take $E=E_0+\kappa_{0.01},\eta=\eta_{0.01}$, 
with Theorem~\ref{thm:counting} and Lemma~\ref{lem:two-H-func}, we have
\begin{align}
    &|H(E,\eta,E_0)\delta\mathfrak n(E_0,t_0)-H(E,\eta,0)\delta\mathfrak n(0,t_0)|<8\frac{\eta_{0.01}^3}{\kappa_{0.01}^3}\times C\frac{\log N}{N}\,,\\
    &|\int_0^{E_0}(\partial_\sigma H(E,\eta,\sigma))\delta\mathfrak n(\sigma,t_0)\mathrm d\sigma|<C\frac{\log N}{N}\int_0^{E_0}\partial_\sigma H(E,\eta,\sigma)\mathrm d\sigma\nonumber\\
    &<8\frac{\eta_{0.01}^3}{\kappa_{0.01}^3}\times C\frac{\log N}{N}\,,
\end{align}
with the scaling $\kappa_{0.01}\sim N^{-\frac23}(\log N)^2$, $\eta_{0.01}\sim N^{-\frac23}(\log N)^\frac32$, we obtain
\begin{align}
    |\int_0^{E_0} H(E,\eta,\sigma)\mathrm d(\mu_N(\sigma,t_0)-\mu_\infty(\sigma,t_0))|<C\frac{1}{N(\log N)^\frac12}\,.
\end{align}

Now we claim that $\sigma_N(t_0)<E_0+2\kappa_{0.01}$ and prove it by contradiction. Suppose that $\sigma_N(t_0)\geq E_0+2\kappa_{0.01}$, the integral
\begin{align}
    \int_{E_0}^\infty H(E,\eta,\sigma)\mathrm d(\mu_N(\sigma,t_0)-\mu_\infty(\sigma,t_0))\geq\frac1N H(E,\eta,\sigma_N(t_0))\,,
\end{align}
where we used the fact that $\mu_\infty(\sigma,t_0)=0$ for $\sigma>E_0$, and that $H(E,\eta,\sigma)$ is non-negative.

By Lemma~\ref{lem:two-H-func}, we have a lower bound
\begin{align}
    H(E_0+\kappa_{0.01},\eta_{0.01},\sigma_N(t_0))&>\pi-4(\frac{\eta_{0.01}}{\sigma_N(t_0)-E_0-\kappa_{0.01}})^3\nonumber\\
    &>\pi-\frac{C}{(\log N)^{\frac32}}>\frac{3\pi}4,
\end{align}
implying that
\begin{align}
    \int_{E_0}^\infty H(E,\eta,\sigma)\mathrm d(\mu_N(\sigma,t_0)-\mu_\infty(\sigma,t_0))\geq\frac{3\pi}{4N}\,,
\end{align}
and
\begin{align}
    |D(E_0+\kappa_{0.01},\eta_{0.01})|\geq \frac{3\pi}{4N}-C\frac1{N(\log N)^\frac12}>\frac{\pi}{2N}\,.
\end{align}
This contradicts with the bound 
\begin{align}
    |D(E_0+\kappa_{0.01},\eta_{0.01})|<\frac{0.03}N
\end{align}
from Corollary~\ref{coro:Delta_bound_outedge}.

\end{proof}

\section{Appendix}\label{sec:appendix}
\subsection{Proof of Lemma~\ref{lem:char-inverse-exist}}
\label{pf:char-inverse-exist}
\begin{proof}
    Consider the inverse ODE: fix $z,t$, let $\Xi(z,s)$ be the solution to
    \begin{align}
        \frac{\mathrm d}{\mathrm ds} \Xi(z,s)=-m_\mu(\Xi(z,s),t-s)\,,&&\Xi(z,0)=z\,,
    \end{align}
    we only need to prove that $\Xi(z,s)$ will stay in $\mathcal S^{++}$ for all $s\in[0,t]$.
    
    First consider the imaginary part, for $z=E+\bm i\eta$,
    \begin{align}
        \Im(-m_\mu(z,t))=&-\Im\Big(\int\frac{\sinh z}{\cosh z-\cosh\sigma}\mathrm d\mu_t(\sigma)\Big)\nonumber\\
        =&\frac12\int\Big(\frac{\sin\eta}{\cosh(E+\sigma)-\cos\eta}+\frac{\sin\eta}{\cosh(E-\sigma)-\cos\eta}\Big)\mathrm d\mu_t(\sigma)\,,
    \end{align}
    $\Im(\Xi(z,s))$ is clearly increasing, and by 
    \begin{align}
        \frac{\sin \eta}{\cosh(E\pm\sigma)-\cos\eta}\leq\frac{\sin\eta}{1-\cos\eta}=\cot\frac\eta2\,,
    \end{align}
    we have $\Im(-m_\mu(z,t))\leq \cot\frac\eta 2$.

    Let $\eta^c_s=\pi-\Im(\Xi(z,s))$, 
    then
    \begin{align}
        \dot \eta^c_s=-|\Im(m_\mu(z,s))|>-\cot\frac{\Im(\Xi(z,s))}2=-\tan\frac{\eta^c_s}2\,,
    \end{align}
    so $\eta^c_s>0$ holds for all $s\in[0,t]$, implying that $\Xi(z,s)<\pi$.

    For the real part,
    \begin{align}
        \Re(m_\mu(z,t))=&\frac12\int(\frac{\sinh(E+\sigma)}{\cosh(E+\sigma)-\cos\eta}+\frac{\sinh(E-\sigma)}{\cosh(E-\sigma)-\cos \eta})\mathrm d\mu_t(\sigma)\nonumber\\
        =&\sinh E\int\frac{\cosh E-\cosh \sigma\cos \eta}{(\cosh(E+\sigma)-\cos \eta)(\cosh(E-\sigma)-\cos \eta)}\mathrm d\mu_t(\sigma)\,,
    \end{align}
    which is proportional to $\sinh E$ times a finite integral, so 
    $\Re(\Xi(z,s))>0$ holds for all $s\in[0,t]$.

    Therefore, the inverse $\Xi(z,t)$ is well-defined. By the fact that $-m_\mu(z,t)$ is holomorphic in $z\in\mathcal S^{++}$ and continuous in $t$, we know that $\Xi(z,t)$ is also holomorphic in $z\in\mathcal S^{++}$.
\end{proof}

\subsection{Proof of Lemma~\ref{lem:m_inf-sol}}
\label{pf:m_inf-sol}
\begin{proof}[Proof of Lemma~\ref{lem:m_inf-sol}]
    We prove that the equation $z=\xi+tm_\infty(\xi,0)$ of $\xi$ has a unique solution $\Xi(z,t)\in\mathcal S^{++}$. Then one can verify the solutions $S_\infty(z,t)$, $m_\infty(z,t)$ and $Z_\infty(\xi,t)$ by implicit differentiation. 
    
    Let $f=m_\infty(\cdot,0)$. We prove the existence of the solution $\xi$ by computing the winding number of the function $\phi_z(\xi)=-z+\xi+tf(\xi)$ for a region $D=[0,R]\otimes[\frac12\Im(z),\pi]$, where $R=2(|z|+t\cot\frac{\Im(z)}{4})$. Notice that $\Re(f(\xi))=0$ if $\Re(\xi)=0$, and $\Im(f(\xi))=0$ if $\Im(\xi)=\pi$.
    
    On $\partial D$, the arguments of $\phi_z(x+\bm iy)$ are:
    \begin{enumerate}
        \item $x\in[0,R],y=\frac12\Im(z)$:
        \begin{align}
            \Im(\phi_z(\xi))=-\Im(z)+\frac12\Im(z)+t\Im(f(\xi))<-\frac12\Im(z)<0\,,
        \end{align}
        \item $x=R,y\in[\frac12\Im(z),\pi]$:
        \begin{align}
            \Re(\phi_z(\xi))=&-\Re(z)+R+t \Re(f(\xi))\nonumber\\
            >&-\Re(z)+R-t \cot\frac y2\nonumber\\
            \geq &-|z|+R-t\cot\frac{\Im(z)}4=\frac R2>0\,,
        \end{align}
        \item $x\in[0,R],y=\pi$: 
        \begin{align}
            \Im(\phi_z(\xi))=-\Im(z)+\pi+0>0\,,
        \end{align}
        \item $x=0,y\in[\frac12\Im(z),\pi]$: 
        \begin{align}
            \Re(\phi_z(\xi))=&-\Re(z)<0\,.
        \end{align}
    \end{enumerate}
    so the winding number is $1$, and there exists a unique solution $\xi$ in the interior of $D$ which is a subset of $\mathcal S^{++}$.

    % The claim that $\xi$ is holomorphic in $z$ is proved by that the derivative $
    % (\frac{\partial \xi}{\partial z})=\frac1{1+tm'(\xi,0)}\neq 0$.

    % In fact, consider $\Re(m'(\xi,0))$, by Lemma~\ref{lem:mz-control-more} and Cauchy-Riemann equation for analytic functions,
    % \begin{align}
    %     \Re(m'(\xi,0))>\frac{\Im(m(\xi,0))}{\Im(\xi)}\,,
    % \end{align}
    % which implies that
    % \begin{align}
    %     \Re(1+tm'(\xi,0))>1+t\frac{\Im(m(\xi,0))}{\Im(\xi)}=\frac{\Im(z)}{\Im(\xi)}>0\,,
    % \end{align}
    % so $\frac1{1+tm'(\xi,0)}\neq 0$.

\end{proof}

\subsection{Proof of Lemma~\ref{lem:im-control}}
\label{pf:im-control}
\begin{proof}
    We first compute $\Theta$. If we replace $m_N(z,t)$ by $m_\infty(z,t)$ in the expression of $\langle \mathrm d\Im(L_N(z,t))\rangle$ in Lemma~\ref{lem:control-martingale}, each term has a closed form when integrated along the characteristic $Z_s$. From the computation of $I_1$, we know
    \begin{align}
        &\int_0^t\partial_zm_\infty(Z_s,s)\mathrm ds=\log(1+tm'(\xi,0))\,,\\
        &\int_0^t\coth(Z_s)m_\infty(Z_s,s)\mathrm ds=\log(\frac{\sinh z}{\sinh \xi})\,,
    \end{align}
    while for the last term, consider $Y_s=\cosh(Z_s)$, it satisfies
    \begin{align}
        \frac{\mathrm d}{\mathrm ds}Y_s=\sinh(Z_s)\frac{\mathrm d}{\mathrm ds}Z_s=\sinh(Z_s)m_\infty(Z_s,s)\,,
    \end{align}
    so
    \begin{align}
        &\int_0^t-\frac{\Im(\sinh(Z_s)m_\infty(Z_s,s))}{\Im(\cosh (Z_s))}\mathrm ds=-\log\frac{\Im(Y_t)}{\Im(Y_0)}=\log\frac{\Im(\cosh \xi)}{\Im(\cosh z)}\,.
    \end{align}

    Then we have
    \begin{align}
        \Big\langle\Im\big(\int_0^t\mathrm dL_N(Z_s,s)\big)\Big\rangle-\frac1\beta\Theta=&\frac1\beta\int_0^t\Big(\Re\big[\partial_z^2\Delta_N(z,s)+\coth(z)\partial_z\Delta_N(z,s)\big]\nonumber\\
        &-\frac{\Im(\sinh(z)\partial_z\Delta_N(z,s))}{\Im(\cosh(z))}\Big)\mathrm ds\,.
    \end{align}

    Now for the integral of the first two terms, using Lemma~\ref{lem:int-eta-2} similar to $I_2$ and $I_3$, we obtain bounds of the form
    \begin{align}
    |\int_0^t\Re\big[\partial_z^2\Delta_N(z,s)
    +\coth(z)\partial_z\Delta_N(z,s)\big]\mathrm ds|<C\frac{(\log N)t}{N\eta\eta_0}\,.
    \end{align}

    For the integral of the last term, use the identity
    \begin{align}
        \frac{\Im(\sinh(z)\partial_z\Delta_N(z,t))}{\Im(\cosh(z))}=\cot\eta\,\Im(\delta_N(z,t))+\coth E\,\Re(\delta_N(z,t))\,,
    \end{align}
    where $\delta_N=\partial_z\Delta_N$, $z=E+\bm i\eta\in\mathcal S^{++}$.

    Since $\cot\eta<\frac 1\eta$ holds for $\eta\in(0,\pi)$, the first integral is also bounded by $C\frac{(\log N)^2t}{N\eta\eta_0}$.

    For the term $\coth E\, \Re(\delta_N(z,t))$, the divergence of $\coth E$ when $E$ approaches $0$ makes a simple bound not available.

    In the case where $E\geq\eta$, use that $\coth E<\coth \eta<1+\frac1\eta<\frac{1+\pi}\eta$, combined with Cauchy's formula, we have 
    \begin{align}
        |\coth E\, \Re(\delta_N(z,t))|<\frac{1+\pi}\eta \frac{C\log N}{N\frac\eta2}=C'\frac{\log N}{N\eta^2}\,.
    \end{align}

    If $E<\eta$, we need to treat the divergence of $\coth E$.
    Use the fact that
    \begin{align}
        \Re(\frac{\sinh z}{\cosh z-\cosh\sigma})=\frac{\sinh E(\cosh E-\cosh\sigma\cos\eta)}{(\cosh(E+\sigma)-\cos\eta)(\cosh(E-\sigma)-\cos\eta)}\,,
    \end{align}
    we know that $\Re(\delta_N(\bm i\eta,t))=0$ (make an extension of $\Delta_N(z,t)$ from $z\in\mathcal S^{++}$ to $z\in\mathcal S^+$ which is even in $\Re(z)$), then
    \begin{align}
        &|\coth E\,\Re(\delta_N(E+\bm i\eta,t))|=|\coth E\,\Re\big(\int_{\bm i\eta}^{E+\bm i\eta} \partial_z \delta_N(z,t)\mathrm dz\big)|\nonumber\\
        &\leq E\coth E\sup_{e\in[0,E]}|\partial_z\delta_N(e+\bm i\eta,t)|\,.
    \end{align}

    Now since $E<\eta<\pi$, $E\coth E$ is bounded by a constant. Take a circle centered at $z$ with radius $\frac\eta2$ across the imaginary axis but still in $\mathcal S^+$, Cauchy's formula gives that
    \begin{align}
        &|\coth E\,\Re(\delta_N(E+\bm i\eta,t))|\leq E\coth E\times 2!\frac{C\log N}{N(\frac\eta2)^2}<C'\frac{\log N}{N\eta^2}\,.
    \end{align}

    In conclusion, we have
    \begin{align}
        |\coth E\,\Re(\delta_N(z,t))|<C\frac{\log N}{N\eta^2}\,,
    \end{align}
    so Lemma~\ref{lem:int-eta-2} implies
    \begin{align}
        |\int_0^t-\coth(E_s)\,\Re(\delta_N(Z_s,s))\mathrm ds|<\int_0^tC\frac{\log N}{N\eta_s^2}\mathrm ds=C\frac{(\log N)t}{N\eta\eta_0}\,.
    \end{align}
    
\end{proof}

\subsection{Proof of Lemma~\ref{lem:nevanlinna-boundary-estimates}}
\label{pf:nevanlinna}
\begin{proof}
Fix \(x\in I\) and \(0<y<d(x)/2\), and abbreviate
\[
    d:=d(x),\qquad r:=\frac{y}{d}<\frac12,
    \qquad z:=x+iy.
\]

We first prove
\eqref{eq:nevanlinna-hyperbolic-arg}
and
\eqref{eq:nevanlinna-hyperbolic-log}.
By the Herglotz representation theorem,
\[
    F(z)
    =
    \alpha+\beta z
    +
    \int_{\mathbb R}
    \left(
        \frac{1}{s-z}-\frac{s}{1+s^2}
    \right)\,d\mu(s),
\]
where \(\alpha\in\mathbb R\), \(\beta\ge0\), and \(\mu\) is a
positive Borel measure satisfying
\[
    \int_{\mathbb R}\frac{d\mu(s)}{1+s^2}<\infty.
\]
Since \(F\) extends holomorphically across \(I\) and is real-valued
there, the Stieltjes inversion formula gives
\[
    \operatorname{supp}\mu\cap I=\varnothing.
\]
In particular,
\[
    |s-x|\ge d,
    \qquad s\in\operatorname{supp}\mu.
\]

Differentiating the representation and taking imaginary parts yield
\[
    F'(z)
    =
    \beta+\int_{\mathbb R}\frac{d\mu(s)}{(s-z)^2}
\]
and
\[
    M:=
    \frac{\operatorname{Im}F(z)}{y}
    =
    \beta+
    \int_{\mathbb R}
    \frac{d\mu(s)}{(s-x)^2+y^2}
    >0.
\]
Consequently,
\begin{equation}
    F^h(x,y)
    =
    \frac{1}{M}
    \left(
        \beta+
        \int_{\mathbb R}
        \frac{s-x+iy}{s-x-iy}\,
        \frac{d\mu(s)}{(s-x)^2+y^2}
    \right).
    \label{eq:nevanlinna-convex-average}
\end{equation}
Thus \(F^h(x,y)\) is a convex average of \(1\) and the points
\[
    \frac{s-x+iy}{s-x-iy},
    \qquad s\in\operatorname{supp}\mu.
\]
Each of these points has modulus one and satisfies
\[
    \left|
        \operatorname{Arg}
        \frac{s-x+iy}{s-x-iy}
    \right|
    =
    2\arctan\frac{y}{|s-x|}
    \le 2\arctan r.
\]
Because \(2\arctan r<\pi/2\), the sector
\[
    \{w\in\mathbb C:
      \operatorname{Re}w>0,\ 
      |\operatorname{Arg}w|\le2\arctan r\}
\]
is convex. Hence
\[
    \left|\operatorname{Arg}F^h(x,y)\right|
    \le2\arctan r<2r,
\]
which proves \eqref{eq:nevanlinna-hyperbolic-arg}.

Moreover, \eqref{eq:nevanlinna-convex-average} gives
\[
    |F^h(x,y)|\le1
\]
and
\[
    \operatorname{Re}F^h(x,y)
    \ge\cos(2\arctan r)
    =\frac{1-r^2}{1+r^2}.
\]
Therefore
\[
    \frac{1-r^2}{1+r^2}
    \le |F^h(x,y)|\le1.
\]
In particular, \(F^h(x,y)\ne0\), and
\[
    0\le-\log|F^h(x,y)|
    \le
    \log\frac{1+r^2}{1-r^2}
    \le
    \frac{2r^2}{1-r^2}
    <
    \frac83r^2.
\]
This proves \eqref{eq:nevanlinna-hyperbolic-log}.

It remains to prove \eqref{eq:nevanlinna-ratio-arg}.
For \(c\in\{-1,1\}\), define
\[
    u_c(z):=\operatorname{Arg}(F(z)-c).
\]
Since \(F(\mathbb H)\subset\mathbb H\), the function \(u_c\) is
harmonic on \(\mathbb H\) and satisfies
\[
    0<u_c(z)<\pi.
\]
The assumption \(F(s)>1\) for \(s\in I\) implies that \(u_c\)
extends continuously to \(I\) with
\[
    u_c(s)=0,\qquad s\in I.
\]
The Poisson representation for bounded harmonic functions therefore
gives
\begin{align*}
    u_c(x+iy)
    &\le
    \int_{\mathbb R\setminus(a,b)}
    \frac{y}{(s-x)^2+y^2}\,ds
    \\
    &=
    \arctan\frac{y}{x-a}
    +
    \arctan\frac{y}{b-x}
    \\
    &\le2\arctan r.
\end{align*}
It follows that
\[
    0<
    \operatorname{Arg}(F(z)-1)
    +
    \operatorname{Arg}(F(z)+1)
    \le4\arctan r.
\]

On the other hand, \(a\ge1\) implies
\[
    x-1\ge x-a\ge d,
    \qquad x+1\ge d.
\]
Hence
\begin{align*}
    0<
    \operatorname{Arg}(z-1)+\operatorname{Arg}(z+1)
    &=
    \arctan\frac{y}{x-1}
    +
    \arctan\frac{y}{x+1}
    \\
    &\le2\arctan r.
\end{align*}
Thus the quantity
\[
    \Delta:=
    \operatorname{Arg}(F(z)-1)
    +
    \operatorname{Arg}(F(z)+1)
    -
    \operatorname{Arg}(z-1)
    -
    \operatorname{Arg}(z+1)
\]
belongs to
\[
    [-2\arctan r,\,4\arctan r]
    \subset(-\pi,\pi).
\]
Since
\[
    \frac{F(z)^2-1}{z^2-1}
    =
    \frac{(F(z)-1)(F(z)+1)}{(z-1)(z+1)},
\]
we conclude that
\[
    \operatorname{Arg}
    \frac{F(z)^2-1}{z^2-1}
    =\Delta.
\]
Consequently,
\[
    \left|
        \operatorname{Arg}
        \frac{F(z)^2-1}{z^2-1}
    \right|
    \le4\arctan r<4r,
\]
proving \eqref{eq:nevanlinna-ratio-arg}.
\end{proof}

\subsection{Proof of Lemma~\ref{lem:two-H-func}}
\label{pf:two-H-func}

\begin{proof}
For $x>0$, set
\[
h(x)
=
2\mathrm{Arg}\sinh\frac{x+i\eta}{2}
-
\mathrm{Arg}\sinh\frac{x+2i\eta}{2}.
\]
Direct differentiation gives
\[
h''(x)
=
\sin\eta\,\sinh x
\left[
\frac{1}{(\cosh x-\cos\eta)^2}
-
\frac{\cos\eta}{(\cosh x-\cos2\eta)^2}
\right]>0.
\]
Since $h(x)\to0$ as $x\to\infty$, it follows that $h$ is
strictly decreasing, strictly convex, and positive.

Moreover, twice integrating with respect to the imaginary part yields
\[
\begin{aligned}
h(x)
&=
\frac{\sinh x}{2}
\int_0^\eta\int_0^\eta
\frac{\sin(u+v)}
     {(\cosh x-\cos(u+v))^2}\,du\,dv\\
&\le
\frac{\eta^3\cosh(x/2)}{4\sinh^3(x/2)}
\le
\frac{\eta^3}{16\sinh^3(x/4)\cosh(x/4)}
\le \frac{4\eta^3}{x^3}.
\end{aligned}
\]
Here we used $\sin t\le t$, the double-angle identities,
and $\sinh t\ge t$.

The factorization
\[
\cosh(E+iy)-\cosh\sigma
=
2\sinh\frac{E+\sigma+iy}{2}
 \sinh\frac{E-\sigma+iy}{2}
\]
gives, with principal arguments,
\[
H(E,\eta,\sigma)
=
\begin{cases}
h(E+\sigma)+h(E-\sigma), & 0\le\sigma<E,\\
\pi+h(E+\sigma)-h(\sigma-E), & \sigma>E.
\end{cases}
\]
Thus, for $\sigma>0$ and $\sigma\ne E$, strict convexity of $h$ gives
\[
\partial_\sigma H
=
h'(E+\sigma)-h'(|E-\sigma|)>0.
\]
Together with continuity at $\sigma=E$, this proves the first claim.

Finally, the bound on $h$ gives
\[
0<H(E,\eta,\sigma)
\le \frac{8\eta^3}{(E-\sigma)^3},
\qquad 0\le\sigma<E,
\]
while its positivity and strict decrease give
\[
0<\pi-H(E,\eta,\sigma)
=h(\sigma-E)-h(\sigma+E)
<h(\sigma-E)
\le \frac{4\eta^3}{(\sigma-E)^3},
\qquad \sigma>E.
\]
This proves the remaining claims.
\end{proof}

\bibliographystyle{amsplain}
\bibliography{bib}
\end{document}